\documentclass[12pt]{article}

\usepackage{amsmath}
\usepackage{amssymb}
\usepackage{latexsym}

\newcommand{\anf}{\raisebox{0.2ex}{\scriptsize$\triangleleft$}}
\newcommand{\ang}{\raisebox{0.2ex}{\scriptsize$\triangleright$}}

\newcommand{\mn}{\medskip\noindent}
\newcommand{\bn}{\bigskip\noindent}
\newcommand{\sn}{\smallskip\noindent}

\newcommand{\lti}{\,{\scriptstyle\ltimes}\,} 
\newcommand{\rti}{\,{\scriptstyle\rtimes}\,} %

\newcommand{\D}{{\mathcal{D}}} 
\newcommand{\E}{{\mathcal{E}}} 
\newcommand{\cR}{{\mathcal{R}}} 
\newcommand{\G}{{\mathcal{G}}} 
\newcommand{\Hh}{{\mathcal{H}}}
\newcommand{\cO}{{\mathcal{O}}}
\newcommand{\K}{{\mathcal{K}}}
\newcommand{\U}{{\mathcal{U}}}

\newcommand{\X}{{\mathcal{X}}}
\newcommand{\cX}{{\mathcal{X}}}
\newcommand{\A}{{\mathcal{A}}}
\newcommand{\cC}{{\mathcal{C}}}
\newcommand{\cU}{{\mathcal{U}}}
\newcommand{\F}{{\mathcal{F}}}
\newcommand{\Q}{{\mathcal{Q}}}
\newcommand{\cP}{{\mathcal{P}}}

\newcommand{\Z}{\mathbb{Z}}
\newcommand{\N}{\mathbb{N}}
\newcommand{\R}{\mathbb{R}}

\newcommand{\dR}{\mathbb{R}}

\newcommand{\C}{\mathbb{C}}
\newcommand{\dC}{\mathbb{C}}

\newcommand{\gM}{{\mathfrak{M}}}
\newcommand{\Lin}{{\mathrm{Lin}}}
\newcommand{\id}{{\mathrm{id}}}
\newcommand{\ad}{{\mathrm{ad}}}

\newcommand{\Tr}{{\mathrm{Tr}}}

\begin{document}
\vspace{10.6mm}
\title{}\date{}
\begin{center}
{\bf {\Large{Hilbert Space Representations of Cross Product Algebras}}}
\renewcommand{\thefootnote}
{\fnsymbol{footnote}}\\
\vspace{43.5pt}{\large Konrad Schm\"udgen
\hspace{40pt}Elmar Wagner}

\bn
\small{Fakult\"at f\"ur Mathematik und Informatik\\ Universit\"at Leipzig, 
Augustusplatz 10, 04109 Leipzig, Germany\\ E-mail: schmuedg@mathematik.uni-leipzig.de / wagner@mathematik.uni-leipzig.de}
\end{center}
{\bf Abstract:} Hilbert space representations of cross product $\ast$-algebras of 
the Hopf $\ast$-algebras $\cU_q(gl_2)$ with the coordinate algebras 
$\hat{\cO}(\C^2_q)$ and $\cO(\R^3_q)$ of quantum vector spaces, and of $\cU_q(su_2)$ with the coordinate algebras $\cO(SU_q(2))$ and 
$\cO(S_q^2)$ of the corresponding quantum spheres, are investigated and 
classified. Invariant states on the coordinate $\ast$-algebras are described 
by two variants of the quantum trace.

\setcounter{section}{0}
\section{Introduction}
Let $\cU$ be a Hopf $\ast$-algebra and let $\X$ be a unital right $\cU$-module 
$\ast$-algebra. Then the smash product algebra $\cU \# \X$ [M] is a 
$\ast$-algebra which contains $\cU$ and $\X$ as $\ast$-subalgebras. Following 
the terminology used by operator algebraists we prefer to call the algebra $\cU \# \X$ a cross product algebra and denote it by $\cU\lti \X$. The cross 
product $\ast$-algebra has a natural physical interpretation: 
We think of $\X$ as an algebra of functions on a ``quantum space" on which 
the elements of $\cU$ act as ``generalized differential operators". Then $\cU\lti \X$ can be considered as an algebra of differential operators with 
coefficients in $\X$ and so as a phase space algebra associated with the 
quantum space. Therefore, as is in usual 
quantum mechanics, Hilbert space representations of the phase space 
$\ast$-algebra $\cU\lti \X$ play a crucial role in the study of the 
quantum space.

Suppose there exists a $\cU$-invariant state $h$ 
on $\X$ and let $\pi_h$ be the GNS-representation of $h$. Then there is a  unique closed
$\ast$-representation of $\cU\lti \X$, called Heisenberg 
representation and denoted also by $\pi_h$, such that 
its restriction to $\X$ is the GNS-representation $\pi_h$ and the cyclic vector 
$v_h$ is $\cU$-invariant. 

In this paper we are concerned with cross product algebras of
the Hopf $\ast$-algebras $\cU_q(gl_2)$ acting on the coordinate algebras 
$\hat{\cO}(\C^2_q)$ and $\cO(\R^3_q)$ of quantum vector spaces, and of
$\cU_q(su_2)$ acting on the coordinate algebras $\cO(SU_q(2))$ and 
$\cO(S_q^2)$ of the corresponding quantum unit spheres. (Precise definitions 
are given in Section 3.) The main purpose of this paper is to study 
Hilbert space representations of these cross product $\ast$-algebras and 
to use them to describe invariant states on 
the coordinate $\ast$-algebras.

%

Rrepresentations of $\cU_q(su_2)\lti\cO(\dR^3_q)$ has been constructed first in [F] and then in [CW] by using other methods. We rediscover these representations in Subsection 6.4.

There are two principal ways to find and to classify Hilbert space 
representations of the $\ast$-algebra $\cU\lti \X$. We discuss these methods and 
some of their applications in 
the case $\cU_q(su_2)\lti \cO(SU_q(2))$. 
In the first approach we assume that the restriction $\pi_{\cU}$ of a $\ast$-representation of $\cU\lti\cX$ to the $\ast$-subalgebra $\cU$ is a "well-behaved" representation of $\cU$.
For $\cU_q(su_2)$ we require that $\pi_\U$ can be expressed as a 
direct sum of spin $l$ representations $T_l$ with arbitrary multiplicities. 
Such $\ast$-representations of $\cU_q(su_2)$ are called integrable. 
For the Heisenberg representation $\pi_h$ of $\cU_q(su_2)\lti\cO(SU_q(2))$ the representation $\pi_{\cU}$ of $\cU_q(su_2)$ is integrable and the Peter-Weyl theorem of $\cO(SU_q(2))$ gives an explicit decomposition of $\pi_\cU$ into a direct sum of representations $T_l$ (see e.g. formula (45) on p. 110 in [KS]).
We prove (Theorem 5.5) that any closed $\ast$-representation of $\cU_q(su_2)\lti \cO(SU_q(2))$ 
obtained from an integrable representation $\pi_\cU$ of $\cU_q(su_2)$ is a direct sum of copies of the Heisenberg representation $\pi_h$. In particular, 
$\pi_h$ is the only closed irreducible 
$\ast$-representation of 
$\cU_q(su_2)\lti \cO(SU_q(2))$ such that its restriction to $\cU_q(su_2)$ is 
integrable. The Heisenberg representation $\pi_h$ is used to 
describe the Haar state $h$ of $\cO(SU_q(2))$ as a quantum trace 
(Theorem 5.7). More precisely, if $C_q$ is the  
Casimir element and $K$ is the group-like generator of $\cU_q(su_2)$, then there exists a holomorphic function 
$\zeta(z)$ on $z \in \C$, ${\rm Re}~ z > 1$, such that
$$
h(x)=\zeta(z)^{-1}~{\rm Tr}~ 
{\overline{\pi_h(C_q)^{-z}\pi_h(K^{-2})}}~ {\overline{\pi_h(x)}},~x \in 
\cO(SU_q(2)).
$$
In the second approach we start with a $\ast$-representation of 
$\X$ written in a canonical form and we try to extend it to a 
$\ast$-representation of $\cU\lti \X$. Since elements of $\cU$ act as 
unbounded operators, during this derivation we have to add some regularity assumptions concerning the unbounded operators. In the case of $\cU_q(su_2)\lti \cO(SU_q(2))$ we end up with irreducible 
representations $(I)_{H,\epsilon}$ parameterized by numbers 
$\epsilon \in \{1,-1\}$ and $H\in (q^{1/2},1]$. 
We prove (Theorem 6.2) that (after extending the domain) the 
representation $(I)_{1,1}$ is unitarily equivalent to the 
Heisenberg representation $\pi_h$. We use the representation 
$(I)_{H,\epsilon}$ to express the Haar state $h$ as 
a partial quantum trace (Theorem 6.4). 

Each of two approaches has its advantage: In the first one the generators of $\cU_q(su_2)$ act in the standard form used in representation theory, while in the second one the action of the generators of $\cO(SU_q(2))$ is given in a canonical form. The second way yields a larger class of representations $(I)_{H,\epsilon}$, because the restriction of $(I)_{H,\epsilon}$ to $\cU_q(su_2)$ is not integrable if $H\in(q^{1/2},1)$.

This paper is organized as follows. In Section 2 we collect 
definitions and general facts on cross product $\ast$-algebras. 
In Section 3 we define the cross product algebras 
studied in this paper and characterize them in terms of generators and 
defining relations. Section 4 contains a number of preliminaries on 
representations. In Subsection 4.2 we list the canonical form of 
representations of the coordinate $\ast$-algebras 
$\hat{\cO}(\C^2_q)$, $\cO(\R^3_q)$, $\cO(SU_q(2))$, and 
$\cO(S_q^2)$ as used in Section 6. The main results of the paper are 
contained in Sections 5 and 6. Section 5 is concerned with 
invariant functionals on the four coordinate $\ast$-algebras and with
Heisenberg representations of the cross product algebras. 
While Section 5 deals with the first approach as explained above, 
the second approach is developed in Section 6. 
For the cross product algebras above and another 
cross product algebra related to the 3D-calculus on $SU_q(2)$ the 
representations are listed by explicit 
formulas of the generators. 

Background material on quantum groups 
can be found  in [FRT], [KS], [M], on unbounded representations in [OS], [P], [S] and on the quantum $SU(2)$ in [KS], [VS], [Wo]. Let us collect some 
definitions and notations used in what follows. The 
comultiplication, the counit and the antipode of a Hopf algebra are denoted 
by $\Delta$, $\varepsilon$, and $S$, respectively. For a 
coaction $\varphi$ and the comultiplication $\Delta$ we freely use 
the Sweedler notations $\varphi(x)=x_{(1)}\otimes x_{(2)}$ and 
$\Delta(x)=x_{(1)}\otimes x_{(2)}$. 
If $T$ is an operator on a Hilbert space, we denote by $\D(T)$ the domain, by 
$\sigma(T)$ the spectrum, by ${\overline T}$ the closure and by $T^\ast$ 
the adjoint of $T$. A self-adjoint operator $A$ is called strictly 
positive if $A{\geq} 0$ and ${\rm ker}~A{=} \{0\}$. We write 
$\sigma (A) \sqsubseteq (a,b]$ if $\sigma(A) {\subseteq} [a,b]$ and 
$a$ is not an eigenvalue of $A$. We say that two self-adjoint operators 
strongly commute if their spectral projections mutually commute.
For a vector $\eta$ of a Hilbert space $\Hh_0$, we denote by 
$\eta_n$ the vector of $\Hh=\oplus^\infty_{k=0}\Hh_k$, $\Hh_k:=\Hh_0$, which has $\eta$ as its n-th component and 
zero otherwise and we put $\eta_{-1}:=0$. 
 By a $\ast$-representation of a $\ast$-algebra $\X$ on a dense domain $\D$ 
of a Hilbert space we mean a homomorphism $\pi$ of $\X$ into the algebra 
$L(\D)$ of linear operators mapping $\D$ into itself such that 
$\langle \pi(x)\eta , \xi \rangle = \langle \eta,\pi(x^\ast)\xi \rangle$ 
for $x \in \X, \eta, \xi \in \D$. A $\ast$-representation $\pi$ is called 
closed if $\D$ is the intersection of domains $\D({\overline{\pi(x)}})$, 
$x \in \X$. 

Throughout this paper we suppose that $q,p\in\dC, p\ne 0$, and $q\ne 0, 1,-1$, and we abbreviate 
\begin{equation*}
\lambda := q-q^{-1}, \gamma := (q+q^{-1})^{1/2},
\lambda_n:=(1-q^{2n})^{1/2},[k]_q:=\lambda^{-1}(q^k-q^{-k}).
\end{equation*}
In Sections 5 and 6 we assume that $0<q<1$, $p>0$, and $p\ne 1.$

\section{Cross product algebras: general concepts}
Throughout this section we suppose that $\cU$ is a Hopf algebra with 
invertible antipode and that 
$\X$ is an algebra (without unit in general).

Let $\X$ be a right $\cU$-module algebra, that is, $\X$ is a right 
$\cU$-module with action $\anf$ satisfying
$$
(xy)\anf f=(x\anf f_{(1)})(y\anf f_{(2)}),~ x,y\in\X, f\in\cU~.
$$
Then the vector space $\cU\otimes \X$ is an algebra, called a
{\it right cross product algebra} and denoted by $\cU\lti\X$, 
 with product defined by
\begin{equation}\label{cross0}
(g\otimes x)(f\otimes y)= gf_{(1)}\otimes (x\anf f_{(2)})y,~ 
x,y\in\X, g,f\in\cU~.
\end{equation}
Let $\cU_0$ be a subalgebra and a right coideal (that is, $\Delta(\cU_0)\subseteq \cU_0\otimes\cU$) of the Hopf algebra $\cU$. 
By (\ref{cross0}), the subspace $\cU_0\otimes\X$ of 
$\cU\otimes \X$ is a subalgebra of $\cU\lti\X$ which we denote 
by $\cU_0\lti\X$. The importance of such algebras stems 
from the fact that the unital subalgebra $\U_0$ generated by the quantum tangent space of a left-covariant differential calculus 
on $\A$ is a right coideal of the Hopf dual $\A^\circ$ ([KS], Proposition 14.5). Since $\A$ is a right 
$\A^\circ$-module algebra (with action given by formula (\ref{rightact}) 
below), the 
algebra $\cU_0\lti \A$ is well defined.

If $\X$ has also a unit element, then we can consider $\X$ and $\U_0$ as 
subalgebras of $\U_0\lti\X$ by identifying $f\otimes 1$ with $f$ and 
$1\otimes x$ with $x$. Then $\U_0\lti\X$ is the 
algebra generated by the two subalgebras $\U_0$ and $\X$ with respect to 
the cross commutation relation
\begin{equation}\label{crossxf}
xf=f_{(1)} (x\anf f_{(2)}),~ x\in\X, f\in\U_0~,
\end{equation}
or equivalently,
\begin{equation}\label{crossfx}
fx=(x\anf S^{-1}(f_{(2)}))f_{(1)},~ x\in\X, f\in\U_0~.
\end{equation}

\sn
Inside the algebra $\U\lti\X$ the right action $\anf$ of $\U$ on $\X$ can 
be nicely expressed by the right adjoint action of the Hopf algebra $\U$. 
Recall that for any $\U$-bimodule $M$ the right adjoint action 
$$
ad_{R}(f)m=S(f_{(1)})m f_{(2)},~ f\in\U,m\in M,
$$ 
is a well defined right action of $\U$ on $M$. The subalgebra $\X$ of 
$\U\lti\X$ is obviously a $\U$-bimodule. 
Using (\ref{crossxf}) we compute 
\begin{equation}\label{adright}
ad_R (f) x {=} S (f_{(1)}) x f_{(2)} {=} S (f_{(1)}) f_{(2)}(x\anf f_{(3)})= \varepsilon (f_{(1)})(x \anf f_{(2)})=x\anf f~.
\end{equation}

Let us turn now to $\ast$-structures. Suppose that $\U$ is a Hopf 
$\ast$-algebra and $\X$ is a right $\U$-module $\ast$-algebra. The latter 
means that $\X$ is a right module algebra and a $\ast$-algebra such that 
the right action $\anf$ and the involution $\ast$ satisfy the compatibility 
condition
\begin{equation}\label{modstar}
(x\anf f)^\ast = x^\ast \anf S(f)^\ast,~ x\in\X, f\in\U~.
\end{equation}
{\bf Lemma 2.1} {\it Let $\U_0$ be $\ast$-subalgebra and a right coideal 
of the Hopf $\ast$-algebra $\U$. Then the algebra $\U_0\lti\X$ is a 
$\ast$-algebra with involution given by }
\begin{equation}\label{fxinv}
(f\otimes x)^\ast:=f^\ast_{(1)}\otimes (x^\ast\anf f^\ast_{(2)}), 
f\in\U_0, x\in\X~.
\end{equation} 
{\bf Proof.} In order to prove that this defines an algebra 
involution we use the formulas (\ref{cross0}) and (\ref{modstar}) and compute
\begin{align*}
((g\otimes x)(f\otimes y))^\ast &=(gf_{(1)}\otimes (x \anf f_{(2)})y )^\ast \\
                                &= f^\ast_{(1)} g^\ast_{(1)}\otimes (y^\ast(x \anf f_{(3)})^\ast)\anf (f^\ast_{(2)} g^\ast_{(2)})\\
                                &=f^\ast_{(1)} g^\ast_{(1)}\otimes (y^\ast \anf (f_{(2)}^\ast g^\ast_{(2)}))(x^\ast \anf (S(f_{(4)})^\ast f_{(3)}^\ast g^\ast_{(3)}))\\
                                &= f^\ast_{(1)} g^\ast_{(1)}\otimes ((y^\ast\anf f^\ast_{(2)})\anf g^\ast_{(2)})(x^\ast \anf g^\ast_{(3)})\\                                 &=(f^\ast_{(1)}\otimes y^\ast\anf f^\ast_{(2)})(g^\ast_{(1)}\otimes x^\ast \anf g^\ast_{(2)})\\
&=(f\otimes y)^\ast (g\otimes x)^\ast~,\\
(f\otimes x)^{\ast\ast}&=(f^\ast_{(1)}\otimes x^\ast\anf f^\ast_{(2)})^\ast=f_{(1)} \otimes (x^\ast\anf f^\ast_{(3)})^\ast \anf f_{(2)} \\
&=f_{(1)}\otimes x\anf (S^{-1}(f_{(3)})f_{(2)})=f\otimes x~.\hspace{2cm}\Box
\end{align*}
If $\X$ has a unit, then the involution (\ref{fxinv}) reads $(f\otimes x)^\ast=(1 \otimes x)^\ast (f\otimes 1)^\ast$.

The above definitions and facts have their left handed 
counterparts. Suppose that $\X$ is a 
left module algebra of  a Hopf algebra $\U$ with left action $\ang$. 
Then the vector space $\X\otimes \U$ is an algebra, called a {\it left cross product algebra} and denoted by $\X\rti\U$, with product defined by
$$
(y\otimes f)(x\otimes g)=y(f_{(1)}\ang x)\otimes f_{(2)} g, ~x,y\in\X,f,g\in\U~.
$$
If $\U_0$ is a subalgebra of $\U$ which is a left coideal 
(i.e. $\Delta(\U_0)\subseteq\U\otimes\U_0)$, then the subspace 
$\X\otimes\U_0$ is a subalgebra of $\X\rti\U$ which is denoted 
by $\X\rti\U_0$. If $\X$ has a unit, then $\X\rti\U_0$ can be considered as the algebra generated by the subalgebras $\X$ and $\U_0$ with cross relation
\begin{equation}\label{oma}
fx=(f_{(1)}\ang x)f_{(2)},~x\in\X, f\in\U_0~.
\end{equation}
If $\U_0$ is a $\ast$-subalgebra and a left coideal of the Hopf $\ast$-algebra $\U$, then the algebra $\X\rti\U_0$ is a $\ast$-algebra with involution defined by
\begin{equation}\label{xfinv}
(x\otimes f)^\ast = (f_{(1)}^\ast \ang x^\ast) \otimes f^\ast_{(2)},~x\in\X, f\in\U_0~.
\end{equation}

All $\U$-module algebras $\X$ occuring in this paper are obtained in the following manner: Let $\A$ be a bialgebra and let $\langle \cdot,\cdot\rangle :\U\times \A\rightarrow\C$ be a dual 
pairing of bialgebras $\U$ and $\A$. If $\X$ is a 
left $\A$-comodule algebra with coaction $\varphi:\X\rightarrow\A\otimes\X$, then $\X$ is a right $\U$-module algebra with right action
\begin{equation}\label{rightact}
x\anf f=\langle f,x_{(1)}\rangle  x_{(2)},~ x\in\X, f\in\U~.
\end{equation}
Then the cross relation (\ref{crossxf}) reads 
\begin{equation}\label{cross2}
xf=f_{(1)} \langle f_{(2)}, x_{(1)}\rangle  x_{(2)},~ x\in\X, f\in\U~.
\end{equation}

If $\X$ is a left comodule $\ast$-algebra of a Hopf $\ast$-algebra $\A$ and 
if $\langle \cdot,\cdot\rangle :\U\times \A\rightarrow\C$ is a dual pairing 
of Hopf $\ast$-algebras, then $\X$ is also a right $\U$-module 
$\ast$-algebra with right action (\ref{rightact}) and so Lemma 2.1 
applies.

Similarly, any right $\A$-comodule algebra $\X$ is a left 
$\U$-module algebra with left action
\begin{equation}\label{leftact}
f\ang x = x_{(1)} \langle f,x_{(2)} \rangle .
\end{equation}
In this case the cross relation (\ref{oma}) can be written as
\begin{equation}\label{cross3}
fx = x_{(1)} \langle f_{(1)},x_{(2)} \rangle f_{(2)}~ x\in\X, f\in\U~.
\end{equation}

\section{Cross product algebras of the Hopf ${\bf \ast}$-al\-ge\-bras 
${\bf \U_q(gl_2)}$ and ${\bf \U_q(su_2)}$}
\subsection{${\bf \cO(M_q(2))}$, ${\bf \cO(SU_q(2))}$, ${\bf \U_q(gl_2)}$ 
and ${\bf \U_q(su_2)}$}
The algebra $\cO(M_q(2))$ has four generators 
$a, b, c, d$ with defining relations
\begin{align}\label{urel1}
a b=qb a, ac=qc a, b d=qdb, cd=qdc, bc = cb,ad-da=\lambda bc~.
\end{align}
The element $\D_q := a d-qbc\equiv da-q^{-1} bc$ is the quantum determinant. 
It is well known that $\cO(M_q(2))$ is a bialgebra and that 
its quotient algebra $\cO(SL_q(2))$ by the 
two-sided ideal generated by $\D_q{-}1$ is a Hopf algebra.

Let $\U_q(gl_2)$ be the algebra with generators $E,F,K,L,K^{-1},L^{-1}$ 
and defining relations
$$
KL=LK,~ KK^{-1}=K^{-1}K=K,~LL^{-1}=L^{-1}L=1 ~,
$$
$$
KEK^{-1}=qE,~KFK^{-1}=q^{-1}F,~LE =EL, ~LF=FL~, $$
$$
EF-FE=\lambda^{-1}(K^2-K^{-2}).
$$
The algebra $\U_q(gl_2)$ is a Hopf algebra with structure maps given by 
\begin{align*}
&\Delta (E)=E\otimes K+K^{-1}\otimes E, 
\Delta(F)=F\otimes K+K^{-1}\otimes F~,\\
&\Delta(K)=K\otimes K, \Delta(L)=L\otimes L, 
\varepsilon(K)=\varepsilon(L)=1, \varepsilon(E)=\varepsilon(F)=0,\\
&S(E)=-qE,S(F)=-q^{-1} F, S(K)=K^{-1}, S(L)=L^{-1}~.
\end{align*}
The Hopf algebra $\U_q(sl_2)$ is the subalgebra of $\U_q(gl_2)$ 
generated by $E$, $F$, $K$ and $K^{-1}$.

There exists a dual pairing $\langle\cdot,\cdot\rangle$ of the 
Hopf algebra $\U_q(gl_2)$ and the bialgebra $\cO(M_q(2))$.
It is determined by the values on the generators $K,L,E,F$ and 
$a,b,c,d$, respectively. The non-zero values are 
\begin{align}\label{dualp1}
&\langle K,a\rangle{=}\langle K^{-1},d\rangle{=}q^{-1/2}, 
\langle K,d\rangle {=}\langle K^{-1},a\rangle{=}q^{1/2},\langle E,c\rangle {=}\langle F,b\rangle{=}1.\\
\label{dualp2}
&\langle L,a\rangle=\langle L,d\rangle=p,~ 
\langle L^{-1},a\rangle =\langle L^{-1},d\rangle=p^{-1}~.
\end{align}
Moreover, (\ref{dualp1}) gives a dual pairing of 
the Hopf algebras $\U_q(sl_2)$ and $\cO(SL_q(2))$.

Suppose in addition that $q$ and $p$ are real. Then $\cO(M_q(2))$ is a $\ast$-bialgebra, 
 $\cO(SL_q(2))$ is a Hopf $\ast$-algebra, denoted by $\cO(SU_q(2))$, 
and $\U_q(gl_2)$ is a Hopf $\ast$-algebra, denoted again 
by $\U_q(gl_2)$, with algebra involutions defined by
\begin{equation}\label{alginv}
a^\ast =d, b^\ast=-qc~{\rm and}~E^\ast=F, 
K^\ast=K, L^\ast=L^{-1}.
\end{equation}
The subalgebra $\U_q(sl_2)$ of $\U_q(gl_2)$ is a Hopf $\ast$-algebra denoted  by $\U_q(su_2)$. Note that the dual 
pairing (\ref{dualp1})--(\ref{dualp2}) of 
$\U_q(gl_2)$ and $\cO(M_q(2))$ satisfies
\begin{equation}\label{starcon}
\overline{\langle S(f)^\ast,x\rangle}=\langle f,x^\ast \rangle,~
 f\in \U_q(gl_2),~ x \in \cO(M_q(2)).
\end{equation}

\subsection{The cross product ${\bf \ast}$-algebras ${\bf \U_q(gl_2)\lti\cO(M_q(2))}$ and ${\bf \cO (M_q(2))\rti \U_q(gl_2)}$}
Since $\cO(M_q(2))$ is a left and right 
comodule algebra  with respect to the comultiplication, $\cO(M_q(2))$ is a 
right and left $\U_q(gl_2)$-module algebra with actions 
given by (\ref{rightact}) and (\ref{leftact}), respectively. 
Hence the cross product algebras $\U_q(gl_2)\lti \cO(M_q(2))$ and 
$\cO(M_q(2))\rti \U_q(gl_2)$ are well defined.

By (\ref{dualp1}), (\ref{dualp2}) and (\ref{cross2}), the generators $E,F,K, L$ and $a, b, c, d$ satisfy 
the following cross relations in the 
{\it right cross product algebra} $\U_q(gl_2)\lti\cO(M_q(2)):$
\begin{align}\label{crgl1}
&aE=q^{-1/2} E a, b E=q^{-1/2} E b~,\\
\label{crgl2}
&cE=q^{1/2} E c+K^{-1} a, d E=q^{1/2} E d+K^{-1} b~,\\
\label{crgl3}
&aF=q^{-1/2} F a+K^{-1} c, b F=q^{-1/2} F b+K^{-1}d~,\\
\label{crgl4}
&cF=q^{1/2} F c, d F=q^{1/2} F d~,\\
\label{crgl5}
&aK =q^{-1/2} K a, bK=q^{-1/2} Kb, c K= q^{1/2} K c,
dK= q^{1/2} K d~,\\
\label{crgl6}
&aL=p La,bL =p Lb, cL = p Lc,dL= p Ld~.
\end{align}
The cross relations of the {\it left cross product algebra} 
$\cO(M_q(2))\rti\U_q(gl_2)$ are:
\begin{align*}
&Ea=q^{1/2} aE+ bK, E b=q^{-1/2}b E~,\\
&Ec=q^{1/2} cE+d K,  Ed=q^{-1/2} E d~,\\
&Fa=q^{1/2} aF, Fb=q^{-1/2}b F+aK~,\\
&Fc=q^{1/2} cF, Fd=q^{-1/2} dF+cK ~,\\
&Ka=q^{-1/2} aK, Kb=q^{1/2} bK, Kc =q^{-1/2} cK, Kd= q^{1/2} dK~,\\
&La=paL ,Lb=p bL, Lc=p cL , Ld=p dL~.
\end{align*}
The two cross product algebras are isomorphic as shown by Lemma 3.1 below. 

Suppose that $q$ and $p$ are real. Since the dual pairing 
of $\U_q(gl_2)$ and $\cO(M_q(2))$ satisfies condition 
(\ref{starcon}), $\cO(M_q(2))$ is a right 
$\U_q(gl_2)$-module $\ast$-algebra. Hence, by Lemma 2.1 and its 
left handed counterpart, 
the cross product algebras $\U_q(gl_2)\lti \cO(M_q(2))$ and 
$\cO(M_q(2))\rti\U_q(gl_2)$ are $\ast$-algebras.

Some simple algebraic facts about the above algebras are collected in the 
next lemma. Its proof is straightforward and will be omitted.

\mn
{\bf Lemma 3.1} (i) {\it The quantum determinant 
$\D_q{=}ad{-}qbc$ is a central element of the algebras} 
$\U_q(su_2){\lti}\cO(M_q(2))$ {\it and} $\cO(M_q(2)){\rti} \U_q(su_2)$.\\
(ii) {\it There is an algebra isomorphism} $\theta$ of 
$\U_q(gl_2){\lti} \cO(M_q(2))$ onto \break
$\cO(M_q(2)){\rti}\U_q(gl_2)$ {\it such that}
\begin{align}\label{thea1}
\theta (a)=a, \theta (d)=d, \theta (b)=-qc, \theta(c)=-q^{-1}b~,\\
\label{theta2}
\theta (E)=F, \theta(F)=E, \theta(K)=K^{-1}, \theta(L)=L^{-1}~.
\end{align}
(iii) {\it If $q$ and $p$ are real, then $\theta$ is a $\ast$-isomorphism of the
$\ast$-algebras}\\
$\U_q(gl_2)\lti \cO(M_q(2))$ {\it and} $\cO(M_q(2))\rti \U_q(gl_2)$.

\sn
Another remarkable property of the isomorphism $\theta$ is that its 
inverse $\theta^{-1}$ is given by the same formulas as $\theta$.

\subsection{The cross product ${\bf \ast}$-algebra ${\bf \U_q(gl_2)\lti\hat{\cO}(\C^2_q)}$}
In the rest of this section we suppose that $q, p \in \R$ , $p\ne 0$, 
and $0<q<1$. 

Let us rename the $\ast$-algebra $\cO(M_q(2))$ defined above by 
$\hat{\cO}(\C^2_q)$ and set $z_1:=b$ and $z_2:=d$. 
By restating the definition of $\cO(M_q(2))$ we see that 
$\hat{\cO}(\C^2_q)$ is the $\ast$-algebra with four generators 
$z_1,z_2,z^\ast_1,z^\ast_2$ and defining relations
\begin{align*}
&z_1z_2=qz_2z_1,~ z_1z_2^\ast=q^{-1}z^\ast_2 z_1, z^\ast_1 z^\ast_2=
q^{-1} z^\ast_2 z^\ast_1,~ z_2z^\ast_1=q^{-1} z^\ast_1 z_2~,\\
&z^\ast_1z_1=z_1z^\ast_1,~ z^\ast_2 z_2 - z_2z^\ast_2=(q^{-2}{-}1)z^\ast_1 z_1~.
\end{align*}
The first equation $z_1z_2=qz_2z_1$ is just the defining relation of the 
coordinate algebra $\cO(\C^2_q)$ of the quantum plane. The $\ast$-algebra 
$\hat{\cO}(\C^2_q)$ introduced above is the left handed 
version of the realification $\cO(\C^2_q)_{Re}$ of the algebra 
$\cO(\C^2_q)$ as defined in [KS], p. 391; see also 
Proposition 9.1.5 therein. 

Let $\C[\cR_q]$ be the $\ast$-algebra of polynomials in a hermitian 
generator $\cR_q$. It is easy to check that there exists a unique injective $\ast$-homomorphism 
$\psi:\hat{\cO}(\C^2_q) \to  \cO(SU_q(2))\otimes\C[\cR_q]$ such that
\begin{equation}\label{hompsi}
\psi (z_1) = b \cR_q,~ \psi (z_2) = d \cR_q~.
\end{equation}
Clearly, $\psi$ is not surjective, because $\cR_q$ is not in the 
image of $\psi$. 
We shall consider $\hat{\cO}(\C^2_q)$ as a $\ast$-subalgebra of $\cO(SU_q(2))\otimes\C[\cR_q]$ by
identifying $x \in \hat{\cO}(\C^2_q)$ with $\psi(x)$. Then we have
\begin{equation}\label{unitsu}  
\cR^2_q = z_1^\ast z_1 + z_2^\ast z_2.
\end{equation}
Thus, $\cO(SU_q(2))$ is just the coordinate algebra of the quantum unit 
sphere, which is obtained  by adding to the relations of $\hat{\cO}(\C^2_q)$ the equation 
$\cR_q^2 = 1$.

In terms of the generators $z_1,z_2,z^\ast_1, z^\ast_2$ the cross 
commutation relations of the {\it right cross product $\ast$-algebra} $\U_q(gl_2)\lti\hat{\cO}(\C^2_q)$ are as follows:
\begin{align}\label{eqp1}
&z_1E=q^{-1/2} Ez_1, z_1^\ast E=q^{1/2} E z_1^\ast - q K^{-1}  z^\ast_2~,\\
\label{eqp2}
&z_2E=q^{1/2} Ez_2+K^{-1} z_1, z_2^\ast E=q^{-1/2} E z_2^\ast ~,\\
\label{eqp3}
&z_1F=q^{-1/2} F z_1 + K^{-1} z_2, z^\ast_1 F= q^{1/2} Fz^\ast_1~,\\
\label{eqp4}
&z_2F=q^{1/2} Fz_2, z^\ast_2 F= q^{-1/2} Fz^\ast_2-q^{-1} K^{-1}z^\ast_1~,\\
\label{eqp5}
&z_1K=q^{-1/2}Kz_1, z^\ast_1K=q^{1/2}Kz^\ast_1, z_2K=q^{1/2} Kz_2, 
z^\ast_2K=q^{-1/2} Kz^\ast_2~,\\
\label{eqp6}
&z_1L=p Lz_1, z^\ast_1L=p Lz^\ast_1, z_2L=p Lz_2, z^\ast_2L=p Lz^\ast_2~.\end{align}
For the element $\cR^2_q\in\hat{\cO}(\C^2_q)$ we have
$\cR^2_qK=K \cR^2_q$ and $\cR^2_q L=p^2 L \cR^2_q$.

\subsection{Two cross product ${\bf \ast}$-algebras containing ${\bf \cO(SU_q(2))}$}
Recall that the Hopf $\ast$-algebra $\cO(SU_q(2))$ is the Hopf 
algebra $\cO(SL_q(2))$ with the involution given by 
$a^\ast= d$ and $b^\ast=-qc$. Hence $\cO(SU_q(2))$ is a right 
$\U_q(su_2)$-module $\ast$-algebra with right action (\ref{rightact}). The corresponding {\it right cross product 
$\ast$-algebra} $\U_q(su_2)\lti\cO(SU_q(2))$ has the cross 
relations (\ref{crgl1})--(\ref{crgl5}). 

We also study another cross product $\ast$-algebra of 
$\cO(SU_q(2))$ where $\U_q(su_2)$ is replaced by a smaller $\ast$-subalgebra.
The quantum tangent space of the $3D$-calculus (see [W] or [KS], p. 407) 
on $SU_q(2)$ is spanned by the elements
$$
X_0:= q^{-1/2} FK, X_2:=q^{1/2} EK, X_1:= (1-q^{-2})^{-1} (1-K^4)
$$
of $\U_q(sl_2)$. These elements satisfy the relations 
$X_0^\ast =X_2$, $X_0^\ast =X_1$ and
\begin{align}\label{tan1}
&q^2 X_1 X_0 - q^{-2} X_0 X_1 = (1 + q^2) X_0~,\\
\label{tan2}
&q^2 X_2 X_1 - q^{-2} X_1 X_2 = (1 + q^2) X_2~,\\
\label{tan3}
&q X_2 X_0 - q^{-1} X_0 X_2 = -q^{-1} X_1~.
\end{align}
Let $\U_0$ be the unital subalgebra of $\U_q(sl_2)$ generated by the 
elements $X_0$,$X_2$,$X_1$. Since the $3D$-calculus is a 
$\ast$-calculus on $SU_q(2)$, 
$\U_0$ is a $\ast$-invariant right coideal of $\U_q(su_2)$, so 
the $\ast$-algebra $\U_0\lti\cO(SU_q(2))$ is well defined. The cross 
relations of the generators $X_0,X_2,Y_1:=1{-}(1{-}q^{-2})X_1$ 
and $a,b,c,d$ of the {\it right cross product 
$\ast$-algebra} $\U_0\lti\cO(SU_q(2))$ read as follows:
\begin{align}\label{cr1}
&aX_2=q^{-1} X_2 a, b X_2=q^{-1}X_2 b~,\\
\label{cr2}
&cX_2=q X_2 c+a, d X_2=qX_2 d+ b~,\\
\label{cr3}
&aX_0=q^{-1} X_0 a+ c, bX_0=q^{-1}X_0 b+d~,\\
\label{cr4}
&cX_0=qX_0 c, d X_0=qX_0 d~,\\
\label{cr5}
&aY_1=q^{-2} Y_1 a,bY_1{=}q^{-2}Y_1b, cY_1= q^{2} Y_1 c,dY_1{=}q^{2}Y_1 d~.
\end{align}

\subsection{The cross product ${\bf \ast}$-algebra ${\bf \U_q(gl_2)\lti\cO(\R^3_q)}$}
Let $\cO(\R^3_q)$ be the $\ast$-algebra with generators $x_1,x_2, x_3$, 
defining relations
\begin{align}\label{r3q}
x_1x_2=q^2x_2x_1,~ x_2x_3=q^2x_3x_2,~ x_3x_1-x_1x_3=\lambda x^2_2,
\end{align}
and involution $x^\ast_1=q^{-1} x_3, x^\ast_2=x_2$, see 
[FRT] or [KS], Proposition 9.14(ii). 

This $\ast$-algebra is a left and right $\cO(M_q(2))$-comodule 
$\ast$-algebra with left coaction 
$\phi_L(x_i)=\sum_j v_{ij}(-q)^{i-j}\otimes x_j$ and right coaction 
$\phi_R(x_i)=\sum_j x_j \otimes v_{ji}$, where 
the matrix $v=(v_{ij})$ is 
$$
v=\left( \begin{matrix} a^2 & q^\prime a b &-b^2\\
                          q^\prime a c & da + qb c & -q^\prime b d\\
-c^2 & -q^\prime c d &d^2 \end{matrix}\right)
$$
and $q^\prime=(1+q^{-2})^{1/2}$. Note that  $v$ is just the 
matrix of the spin 1 corepresentation of $\cO(SU_q(2))$ when the 
elements $a,b,c,d$ are taken as generators of $\cO(SU_q(2))$. Thus, 
$\cO(\R^3_q)$ is a right and left $\U_q(gl_2)$-module $\ast$-algebra with 
actions given by (\ref{rightact}) and (\ref{leftact}). 

The {\it right cross product $\ast$-algebra} $\U_q(gl_2)\lti \cO(\R^3_q)$
has the cross relations:
\begin{align}\label{cr3e}
&x_1 E=q^{-1}E x_1, x_2E= Ex_2-q\gamma K^{-1}x_1,
x_3E=q Ex_3+q\gamma K^{-1}x_2,\\
\label{cr3f}
&x_1F= q^{-1}Fx_1-q^{-1}\gamma K^{-1}x_2,x_2F=Fx_2+q^{-1}\gamma K^{-1}x_3,
x_3F=qFx_3,\\
\label{cr3k}
&x_1K =q^{-1}Kx_1,~ x_2K=Kx_2,~ x_3K=qKx_3,~\\
\label{cr3l}
& x_1L=p^2 Lx_1,~ x_2L=p^2Lx_2,~x_3L = p^2Lx_3~.
\end{align}

Some properties of the algebra $\cO(\R^3_q)$ are collected in the 
following lemma. We omit the proof.

\sn
{\bf Lemma 3.2}
(i) $\Q_q^2:=q^{-1} x_1x_3+qx_3x_1+x^2_2=(1+q^{-2})x_3^\ast x_3+q^2 x^2_2$ {\it belongs to the center of the algebra} $\cO(\R^3_q)$ and 
$\phi_L(\Q_q^2)=\D_q^2\otimes \Q_q^2$.\\
(ii) {\it There is an injective $\ast$-homomorphism} 
$\rho$ {\it of} $\cO(\R^3_q)$ {\it into} $\hat{\cO}(\C^2_q)$ 
{\it such that}\break
$\rho(x_1){=}(1{+}q^2)^{-1/2}(z_2^{\ast 2}{-}q^{-1}z_1^2)$,
 $\rho(x_2){=}z^\ast_1z_2^\ast{+}z_2z_1$,
 $\rho(x_3){=}(1{+}q^2)^{-1/2} (qz^2_2{-}z^{\ast 2}_1).$\\
(iii) $(id{\otimes} \rho)  {\circ} \phi_L{=}\Delta{\circ}\rho$, 
{\it where $\Delta$ is the comultiplication of} 
$\cO(M_q(2)) {=}\hat{\cO}(\C^2_q)$.\\
(iv) {\it There is a $\ast$-isomorphism $\vartheta$ of 
$\U_q(gl_2)\lti \cO(\R^3_q)$ onto $\cO(\R^3_q) \rti \U_q(gl_2)$ such that $\vartheta(x_i)=x_i$, $\vartheta(E)=F$, $\vartheta(F)=E$, 
$\vartheta(K)=K^{-1}$ and $\vartheta(L)=L^{-1}$.}
\medskip

From Lemma 3.2,(ii) and (iii), it follows that $\cO(\R^3_q)$ and 
$\U_q(gl_2)\lti \cO(\R^3_q)$ are $\ast$-subalgebras of 
$\hat{\cO}(\C^2_q)$ and $\U_q(gl(2) \lti\hat{\cO}(\C^2_q)$, respectively, 
if we identify $x\in\cO(\R^3_q)$ and $\rho (x)\in \hat{\cO}(\C^2_q)$. 

The quotient $\ast$-algebra $\cO(S^2_q)$ of $\cO(\R^3_q)$ by the  
ideal generated by the central hermitian element $\Q_q^2{-}1$ is 
called the coordinate algebra of the quantum unit sphere $S^2_q$ of $\R^3_q$.
Let $y_i$ denote the image of the generator $x_i$ of $\cO(\R^3_q)$ under 
the quotient map. The defining relations of the algebra $\cO(S^2_q)$ are
\begin{align}
&y_1y_2=q^2y_2y_1,~ y_2y_3=q^2y_3y_2,~ y_3y_1-y_1y_3=\lambda y^2_2,\notag\\
\label{unitspr} 
&q^{-1}y_1y_3 + qy_3y_1 + y_2^2 =1~.
\end{align}
The quantum sphere $S^2_q$ of $\R^3_q$ is one of Podles' 
quantum spheres $S^2_{qc}$. More precisely, it is the 
quantum sphere $S^2_{q\infty}$ in [Po] and $S^2_{q0}$ in [KS], Section 4.5. 

Since $\phi_L(\Q_q^2)=\D_q^2\otimes \Q_q^2$ by Lemma 3.2(i), the left 
coaction $\phi_L$ of $\cO(M_q(2))$ on $\cO(\R^3_q)$ passes to a left coaction 
of the quotient algebras $\cO(SU_q(2))$ on $\cO(S^2_q)$. Thus, 
$\cO(S^2_q)$ is a left $\cO(SU_q(2))$-comodule and hence right 
$\U_q(su_2)$-module $\ast$-algebra. The corresponding 
{\it right cross product $\ast$-algebra}\break $\U_q(su_2)\lti \cO(S^2_q)$
has the cross relations (\ref{cr3e})--(\ref{cr3k}) with $x_i$ replaced by
$y_i$. 

The assertions (ii)--(iv) of Lemma 3.2 have their 
counterparts for the quantum sphere $S^2_q$. 
The map $\vartheta$ from Lemma 3.2(iv) passes to a $\ast$-isomor\-phism of the 
right and left cross product $\ast$-algebras $\U_q(su_2)\lti \cO(S^2_q)$ 
and 
$\cO(S^2_q)\rti \U_q(su_2)$. There exists an injective $\ast$-isomorphism 
$\rho$ from $\cO(S^2_q)$ into $\cO(SU_q(2))$ such that 
$(id\otimes \rho)\circ \phi_L=\Delta\circ\rho$ and 
$$
\rho(y_1)=(1+q^2)^{-1/2}(a^2{-}q^{-1}b^2),\rho(y_2)=db-qca,
\rho(y_3)=(1+q^2)^{-1/2} (qd^2{-}q^2c^2).
$$
Thus, if we identify $x \in \cO(S^2_q)$ with $\rho(x)\in \cO(SU_q(2))$, 
then $\U_q(su_2)\lti \cO(S^2_q)$ becomes a $\ast$-subalgebra of $\U_q(su_2)\lti \cO(SU_q(2))$.
 
Let $\C[\Q_q]$ be the $\ast$-algebra of polynomials in $\Q_q=\Q_q^\ast$. 
There is an injective $\ast$-homomorphism $\psi:\cO(\R^3_q) \to  
\cO(S_q^2)\otimes\C[\Q_q]$ such that
$\psi (x_i) = y_i \Q_q$, $i=1,2,3$. By identifying $x \in \cO(\R^3_q)$ with 
$\psi (x)$,  $\cO(\R^3_q)$ becomes a $\ast$-subalgebra of 
$\cO(S_q^2)\otimes\C[\Q_q]$. 

\section{Preliminaries on representations}
\subsection{Three auxiliary lemmas}
In this subsection $q$ is a positive real number such that $q\neq 1$. 

\mn
{\bf Lemma 4.1} {\it Up to unitary equivalence each isometry $w$ 
on a Hilbert space $\Hh$ is of the following form: There exist Hilbert 
subspaces $\Hh^u$ and $\Hh^s_0$ of $\Hh$ and a unitary operator $w_u$ on 
$\Hh^u$ such that $w = w_u\oplus w_s$ on $\Hh = \Hh^u \oplus \Hh^s$, where 
$\Hh^s =\mathop{\oplus}^{\infty}_{n=0}\Hh^s_n$, $\Hh^s_n = \Hh^s_0$, 
and $w_s \eta_n = \eta_{n+1}$, $\eta_n \in \Hh^s_n$. This decomposition of $w$ 
is unique, because $\Hh^u =\cap^\infty_{n=0} w^n\Hh$ and $\Hh^s$ is the closed 
linear span of $w^n({\rm ker}~w^\ast), n \in \N_0$.}

\sn
Lemma 4.1 is proved in [SF], Theorem 1.1. The decomposition 
$w = w_u\oplus w_s$ is called the {\it Wold decomposition} of $w$. The operator $w_u$ is the unitary part of $w$ and $w_s$ is a unilateral shift operator of multiplicity $\dim (\Hh^s_0)$. 

\mn
{\bf Lemma 4.2} {\it  Let $A$ be a self-adjoint operator 
and let $w$ be an isometry on a Hilbert space $\Hh$ such that
\begin{equation}\label{waw}
qwA \subseteq A w.
\end{equation}
Then the Wold decomposition $w = w_u\oplus w_s$ on $\Hh = \Hh^u \oplus \Hh^s$
reduces the operator $A$, that is, there are self-adjoint operators $A^u$ on 
$\Hh^u$ and $A^s$ on $\Hh^s$ such that $A= A^u \oplus A^s$, and we have:\\
(i) If $A$ is strictly positive, then there exists a self-adjoint operator $A^u_0$ on a Hilbert space $\Hh^u_0$ 
with $\sigma (A^u_0) \sqsubseteq (q,1]$ if $q<1$ and 
$\sigma (A^u_0) \sqsubseteq (q^{-1},1]$ if $q>1$ such that, up to 
unitary equivalence,
$\Hh^u=\mathop{\oplus}^{\infty}_{n=- \infty}\Hh^u_n$, where 
$\Hh^u_n:=\Hh^u_0$,
$$ 
A^u \eta_n = q^{n} A^u_0 \eta_n,~w_u\eta_n = \eta_{n+1},~~
{\rm for}~ \eta_n \in \Hh^u_n, n \in \Z~.
$$
(ii) There is a self-adjoint operator $A^s_0$ on the Hilbert space $\Hh^s_0$ 
such that 
$$
A^s \eta_n = q^{n} A^s_0 \eta_n,~w_s\eta_n = \eta_{n+1},~~
{\rm for}~ \eta_n \in \Hh^s_n, n \in \N_0~.
$$}
{\bf Proof.} From the functional calculus of self-adjoint operators it follows that equation (\ref{waw}) implies that 
$w\varphi(qA) = \varphi(A)w$ for $\varphi\in L^\infty (\R)$.
Therefore, for any $\varphi\in L^\infty (\R)$, the subspace $\Hh^u:=\cap^\infty_{n=0} w^n \Hh$ is $\varphi (A)$-invariant, and therefore so is its orthogonal complement $\Hh^s$.
This implies that $A$ decomposes as $A= A^u \oplus A^s$ with respect to the 
orthogonal sum $\Hh=\Hh^u \oplus \Hh^s$.

Let $e(\mu)$ denote the spectral 
projections of $A^u$. Let $\Hh^u_n:=e((q^{n+1},q^n])\Hh^u$ and 
$A^u_n:= A^u \lceil \Hh^u_n$, $n \in \Z$. 
Since $A^u$ is strictly positive,  
$\Hh^u=\mathop{\oplus}^{\infty}_{n=- \infty}\Hh^u_n$. 
Since $w_u$ is unitary, (\ref{waw}) implies that $A^u= q w_uA^u w_u^\ast$ and hence 
$\varphi(A^u) = w_u\varphi(qA^u)w_u^\ast$ for $\varphi\in L^\infty (\R)$. 
This yields $w_u\Hh^u_n =\Hh^u_{n+1}$. Thus, 
up to unitary equivalence, we can assume that
$\Hh^u_n =\Hh^u_0$ and $w_u\eta_n =\eta_{n+1}$ for $\eta_n \in \Hh^u_n$. 
Thus,  $A^u\eta_n =q^n w_u^nA^uw_u^{\ast n}\eta_n= q^nw_u^nA^u_0\eta_0 
=q^nA^u_0\eta_n$ which proves (i).

Since $w_s\varphi(qA^s) = \varphi(A^s)w_s$ for $\varphi\in L^\infty (\R)$ by 
(\ref{waw}), $\varphi(A^s)$ leaves ${\rm ker}~w_s^{n\ast} =
\Hh^s_0+\cdots +\Hh^s_{n-1}$ invariant. Since this is true for all 
$\varphi\in L^\infty (\R)$, it follows that $\varphi(A^s)$ leaves
each space $\Hh^s_m$ invariant. Setting $A^s_0\:= A^s \lceil \Hh^s_0$, 
relation (\ref{waw}) gives $A^s\eta_n =q^n A^s_0 \eta_n$. This 
proves (ii). 
\hfill $\Box$

\sn
{\bf Lemma 4.3} {\it Let $x$ be a closed operator on Hilbert space $\Hh$.
Then we have $\D(xx^\ast)=\D(x^\ast x)$, this domain is dense in $\Hh$, 
and the relation 
\begin{equation}\label{qhyp}
x x^\ast -q^2 x^\ast x = 1~
\end{equation} 
holds if and only if $x$ is unitarily equivalent to an orthogonal direct sum of operators of the following forms:}
\begin{align*}
&(I) ~for~any~ q> 0: x\eta_n = ((1-q^{2n})/(1-q^2))^{1/2} \eta_{n-1}~\\ 
&on~\Hh={\oplus}^{\infty}_{n=0}\Hh_n,~\Hh_n =\Hh_0.\\
&(II)_A ~ for~ 0{<}q{<}1: x\eta_n =(1-q^2)^{-1/2}\alpha_{n}(A) \eta_{n-1}\\
&on~\Hh{=}{\oplus}^{\infty}_{n=-\infty}\Hh_n, \Hh_n{=}\Hh_0,
~where~A~is~a~{self\mbox{-}adjoint}~operator~on~the\\ 
&Hilbert~space~\Hh_0 ~such~ that~ \sigma(A) \sqsubseteq (q^2,1]~ and~ \alpha_{n}(A):= (1+q^{2n} A)^{1/2}.\\
&(III)_u~for~0{<}q{<}1:x{=}(1{-}q^2)^{-1/2}u,~where~u~is~a~unitary~operator~on~\Hh.
\end{align*}
{\bf Proof.} Clearly, the operators $x$ have the stated properties. That all 
{\it irreducible} operators $x$ are of one of the above form was proved in [CGP]. 
The general case follows by decomposition theory or by 
modifying the proof in [CGP]. \hfill $\Box$

\sn
Representation $(I)$ is usually called the {\it Fock representation}. 
Note that for $q>1$ the Fock representation is the only representation of 
relation (\ref{qhyp}).  
 
\subsection{Representations of the coordinate ${\bf \ast}$-algebras} 

Suppose that $0<q<1$.
Let $\Hh_0$ and $\G$ be Hilbert spaces and set 
$\Hh=\mathop{\oplus}^{\infty}_{n=0}\Hh_n$, 
where $\Hh_n:=\Hh_0$. 
One checks that the following formulas define $\ast$-representations of the  
corresponding $\ast$-algebras on the Hilbert space $\G\oplus\Hh$: 
\begin{align}
\cO(SU_q(2)):~&a=v,~d=v^\ast,~b=c=0~~{\rm on }~\G~,\hspace{6cm}\nonumber\\
\label{opsu2}
               &a\eta_n{=}\lambda_n\eta_{n-1},
	       d\eta_n{=}\lambda_{n+1} \eta_{n+1}, 
	       b\eta_n{=}q^{n+1} w\eta_n,c\eta_n\makebox[0pt][l]{${=}{-q^n}w^\ast\eta_n,$}
\end{align}
where $w$ and $v$ are unitary operators on 
$\Hh_0$ and $\G$, respectively.
\begin{align}
{\hat\cO}(\C^2_q):~ z_1 =&\  z^\ast_1=0,~z_2=M,~z^\ast_2=M^\ast~~{\rm on}~\G~,\nonumber\\
\label{opc2}
                   z_1\eta_n{=}&\,q^{n{+}1} wA\eta_n,
		   z^\ast_1\eta_n{=}q^{n{+}1}w^\ast A\eta_n, 
		   z_2\eta_n{=}\lambda_{n{+}1} A\eta_{n{+}1},
		   z^\ast_2\eta_n{=}\lambda_nA\eta_{n{-}1},
\end{align}
where $M$ is a normal operator on $\G$, $A$ is a strictly positive 
self-adjoint operator and $w$ is a unitary on $\Hh_0$ such that 
$wAw^\ast=A$.
\begin{align}
\cO(S^2_q):~ & y_1= (1+q^2)^{-1/2} u^\ast,~y_2=0,~ 
y_3=(1+q^{-2})^{-1/2} u~~{\rm on }~\G~,\hspace{2cm}\nonumber\hspace{1cm}\\
\label{opr3}
              & y_1\eta_n=(1+q^2)^{-1/2} \lambda_{2n} \eta_{n-1}, 
y_2\eta_{n}=q^{2n+1} w\eta_n~,\hspace{1cm}\\
\label{opr4}
              & y_3\eta_n=q(1+q^2)^{-1/2} \lambda_{2(n+1)} \eta_{n+1}~,
\end{align}
where $u$ is a unitary operator on $\G$ and $w$ is a self-adjoint unitary 
on $\Hh_0$. \hspace{1cm}
\begin{align}
\cO(\R^3_q):~~ & x_1=q^{-1} M^\ast,~x_2=0,~ 
x_3=M~~{\rm on }~\G~,\hspace{4.7cm}\nonumber\\
\label{opr3}
              & x_1\eta_n=(1+q^2)^{-1/2} \lambda_{2n} A\eta_{n-1}, 
x_2\eta_{n}=q^{2n+1} wA\eta_n~,\\
\label{opr4}
              & x_3\eta_n=q(1+q^2)^{-1/2} \lambda_{2(n+1)} A\eta_{n+1}~,
\end{align}
where $M$ is a normal operator on $\G$ and $A$ is a strictly positive 
self-adjoint operator and $w$ is a self-adjoint unitary on $\Hh_0$ such that 
$wAw^\ast=A$.

\mn
{\bf Lemma 4.4} {\it Any $\ast$-representation of $\cO(SU_q(2))$ or 
$\cO(S_q^2)$ is up to unitary equivalence of the above form.}

\sn 
{\bf Proof.} For $\cO(SU_q(2))$ the assertion is Proposition 4.19 in [KS]. 
We sketch the proof for $\cO(S_q^2)$. The third defining relation of 
$\cO(S_q^2)$ yields
\begin{equation}\label{y3y2}
y_3^\ast y_3 - y_3y_3^\ast = (1-q^2)y_2^2~.
\end{equation}
Combining the latter with (\ref{unitspr}), we obtain
\begin{equation}\label{y3y3}
y_3^\ast y_3 - q^4 y_3 y_3^\ast = q^2(1-q^2)~.
\end{equation}
By (\ref{unitspr}), $y_3$ and $y_2$ are bounded. Hence $\G:={\rm ker}~y_2$ is reducing and $u:= (1+q^{-2})^{1/2}y_3 \lceil  \G$ is unitary on $\G$ by (\ref{y3y2}) and (\ref{y3y3}). 

Assume now that ${\rm ker}~y_2= \{0\}$. The 
representations of relation (\ref{y3y3}) are given by Lemma 4.3. 
For the series $(II)$ and $(III)$ the operator $y_3^\ast y_3 - y_3y_3^\ast$ 
is not strictly positive which contradicts (\ref{y3y2}). Thus, only the Fock 
representation $(I)$ is possible, so $y_3$ has the above form. Let 
$y_2 =v|y_2|$ be the polar decomposition of $y_2$. Since
$y_2$ is self-adjoint and ${\rm ker}~y_2 =\{0\}$, $v$ is self-adjoint and
unitary. From the formula for $y_3$ and (\ref{y3y2}) we obtain 
$y_2^2\eta_n =q^{4n+2}\eta_n$ and so $|y_2|\eta_n=q^{2n+1}\eta_n$. Since 
$v|y_2|v^\ast=|y_2|$, $v$ leaves each space $\Hh_n$ invariant. Hence there are self-adjoint unitaries $v_n$ on $\Hh_n$ such that $v\eta_n=v_n\eta_n$. From the relation $y_2y_3=q^2y_3y_2$ it follows that $w:=v_0=v_n$ for all $n$. \hfill $\Box$   

\sn 
Let us turn now to the $\ast$-algebras ${\hat\cO}(\C^2_q)$ and $\cO(\R^3_q)$. 
Suppose we have a $\ast$-representation of $\cO(SU_q(2))$ resp.\ $\cO(S_q^2)$ 
and a (possibly unbounded) self-adjoint operator $R$ commuting 
with all representation operators. Obviously, we then obtain a 
$\ast$-representation of $\cO(SU_q(2))\otimes\C[\cR_q]$ resp.\ 
$\cO(S_q^2)\otimes\C[\Q_q]$ and so of its $\ast$-subalgebra 
$\hat{\cO}(\C^2_q)$ resp.\ $\cO(\R^3_q)$ on the domain 
$\D:=\cap^\infty_{n=0} \D(R^n)$ which maps $\cR_q$ resp.\ $\Q_q$ into 
$R \lceil \D$. We shall think of $\cR_q$ and $\Q_q$ as 
{\it quantum radii} of the quantum vector spaces $\C^2_q$ and $\R^3_q$, 
respectively. Thus it is natural to require that these elements 
are represented by a positive operator $R$. Let us call 
$\ast$-representations of $\hat{\cO}(\C^2_q)$ resp.\ $\cO(\R^3_q)$ of this form {\it admissible}. 

It is easy to show that the $\ast$-representations of $\hat{\cO}(\C^2_q)$ 
and $\cO(\R^3_q)$ listed above are 
admissible and that any admissible $\ast$-representation of 
$\hat{\cO}(\C^2_q)$ resp.\ $\cO(\R^3_q)$ is up to unitary equivalence of the above form. 

\subsection{Integrable representations of ${\bf \U_q(su_2)}$}
For later use let us restate some well known facts 
(see e.g.\ [KS]). The irreducible unitary corepresentations of 
$\cO(SU_q(2))$ are labeled by numbers $l\in\frac{1}{2}\N_0$. 
By (\ref{leftact}), each such corepresentation gives an irreducible 
$\ast$-representation $T_l$ of $\U_q(su_2)$.  The $\ast$-representation 
$T_l$ acts on a $(2l+1)$-dimensional Hilbert space $V_l$ with orthonormal 
basis $\{e_j;j=-l,-l+1,{\dots},l\}$ by the formulas:
\begin{equation}\label{tlkef}
T_l(K) e_j=q^je_j, T_l(E) e_j=\alpha_{j+1,l} e_{j+1}, T_l(F) e_j=\alpha_{j,l} e_{j-1},
\end{equation}
where $e_{l+1}{=} e_{-l-1}{=}0$ and 
$\alpha_{jl}:=([l{+}j]_q[l-j+1]_q)^{1/2}$. For the Casimir element
\begin{equation}\label{casi}
C_q:=EF+\lambda^{-2}(q^{-1}K^2+qK^{-2}-2)=FE+\lambda^{-2}(qK^2+q^{-1}K^{-2}-2)
\end{equation}
we have
\begin{equation}\label{tlc}
T_l(C_q)=[l+1/2]^2_q=\lambda^{-2}(q^{l+1/2}-q^{-l-1/2})^2.
\end{equation}
In representation theory (see [Wa], Ch. 5 or [S], Ch. 10), 
a $\ast$-representation of a universal enveloping algebra is called 
integrable if it comes from a unitary representation of 
the corresponding connected simply connected Lie group. This suggest the 
following definition.

We say that a closed $\ast$-representation of $\U_q(su_2)$ on a Hilbert 
space is {\it integrable} if it is a direct sum of
 $\ast$-representations $T_l,l\in{\frac{1}{2}}\N_0$. It can be 
shown that a closed $\ast$-representation $\pi$ of $\U_q(su_2)$ is 
integrable if and only if $\pi$ is a direct sum of finite dimensional 
$\ast$-representations and $\sigma(\overline{\pi(K)})\subseteq[0,+\infty)$.

\section{Invariant functionals and Heisenberg representations}

In the rest of this paper we assume that $0<q<1$, $p>0$ and $p \ne1$.

\subsection{Invariant functionals on coordinate algebras}
Let $\X$ be a right module algebra of a Hopf algebra $\U$ with right 
action $\anf$. A linear functional $h$ on $\X$ is called 
$\U$-{\it invariant} if $h(x\anf f)=\varepsilon (f)h(x)$ for 
all $x\in\X$ and $f\in\U$.

Suppose $\X$ is a left comodule algebra of a Hopf algebra $\A$.
A linear functional $h$ on 
$\X$ is said to be $\A$-{\it invariant} if $(\id\otimes h)\varphi(x)=h(x)1$ 
for $x\in\X$ or equivalently if $h(x_{(2)})x_{(1)}=h(x)1$ for $x\in\X$. 

If $\langle\cdot,\cdot\rangle:\U\times \A\rightarrow\C$ is a dual pairing of 
Hopf algebras, then $\X$ is a right $\U$-module algebra with right 
action (\ref{rightact}). In this case the two invariance concepts are not 
equivalent, but they are related as follows.

\mn
{\bf Lemma 5.1} {\it Let $h$ be a linear functional on $\X$. \\
(i) If $h$ is $\A$-invariant, then $h$ is also $\U$-invariant.\\
(ii) Suppose that $\U$ separates the points of $\A$, that is, if $x\in\A$ 
and $\langle f,x\rangle=0$ for all $f\in\U$, then $x=0$. Then $h$ is 
$\A$-invariant if $h$ is $\U$-invariant.}

\mn
{\bf Proof.} (i) is obvious. We verify (ii) and assume that $h$ is 
$\U$-invariant. By the $\U$-invariance and (\ref{rightact}), we have 
$$
h(x\anf f)=\langle f,h(x_{(2)})x_{(1)}\rangle=\varepsilon (f)h(x)=
\langle f,h(x) 1\rangle
$$
for all $f\in\U$ and so $h(x_{(2)})x_{(1)}=h(x)1$.\hfill $\Box$

\mn
Now we specialize the preceding to the case $\U=\U_q(su_2)$ and 
$\A=\cO(SU_q(2))$ with the dual pairing given by (\ref{dualp1}). It is well known 
(see [KS], p. 113 and p. 128) that there is a unique $\A$-invariant 
linear functional $h$ such that $h(1)=1$ on each of the right $\A$-comodule algebras 
$\X=\cO (SU_q(2))$ and $\X=\cO(S^2_q)$. This functional $h$ is a state on the $\ast$-algebra $\X$. For 
$\X=\cO(SU_q(2))$, it is called the Haar state of $SU_q(2)$ and  
explicitly given by
\begin{equation}\label{haar}
h(a^r b^k c^l) = h (d^r b^k c^l) = 
\delta_{r0} \delta_{kl} (-1)^k [k+1]^{-1}_q, r, k, l \in \N_0.
\end{equation} 
Since $\U_q(su_2)$ separates the points of 
$\cO(SU_q(2))$ (see [KS], 4.4), Lemma 5.1 applies and by the 
preceding we have

\mn
{\bf Lemma 5.2} {\it There is a unique $\U_q(su_2)$-invariant linear 
functional $h$ on $\X=\cO(SU_q(2))$ resp.\ $\X=\cO(S^2_q)$ satisfying $h(1)=1$.}\\

Next we look for $\U_q(gl_2)$-invariant linear functionals on 
$\hat{\cO} (\C^2_q)$ and $\cO(\R^3_q)$. We carry out the construction 
for $\hat{\cO}(\C^2_q)$. Replacing $\cO(SU_q(2))$ by $\cO(S^2_q)$ and 
$\cR^2_q$ by $\Q^2_q$, the case of $\cO(\R^3_q)$ is treated completely 
similarly.  Recall that $\hat{\cO}(\C^2_q)$ is a $\ast$-subalgebra of 
$\cO(SU_q(2))\otimes \C[\cR_q]$ via the $\ast$-iso\-mor\-phism $\psi$ 
defined by (\ref{hompsi}). In order 
to ``integrate over the quantum plane", we need more functions of 
$\cR_q$ than polynomials. Let $\F$ denote the $\ast$-algebra of all 
Borel functions on $(0,+\infty)$. For $\varphi\in\F$ we write 
$\varphi(\cR_q)$ instead of $\varphi(t)$ and consider $\C[\cR_q]$ as a
$\ast$-subalgebra of $\F$. Then, $\hat{\cO}_e(\C^2_q):=\cO(SU_q(2))\otimes\F$ is a 
right $\U_q(gl_2)$-module $\ast$-algebra under the actions
$(\varphi\anf L)(\cR_q)=\varphi(p\cR_q)$ and 
$\varphi\anf f=\varepsilon (f)\varphi$ 
for all $\varphi\in\F$ and $f\in\U_q(su_2)$. The cross 
product $\ast$-al\-ge\-bra $\U_q(gl_2)\lti \hat{\cO}_e(\C^2_q)$ has the cross 
relation $\varphi (\cR_q)L=L\varphi(p\cR_q)$ and contains 
$\U_q(gl_2)\lti\hat{\cO}(\C^2_q)$ as $\ast$-subalgebra.

 Let $\mu_0$ be a positive Borel measure supported 
on the interval $(p,1]$ if $p<1$ resp.\ $(p^{-1},1]$ if $p>1$ . Let $C_c(0,+\infty)$ denote the continuous 
function on  $(0,+\infty)$ with compact support. There is a unique positive Borel measure 
$\mu$ on $(0,+\infty)$ such that $\mu(p\gM)=p\mu(\gM)$ for any 
Borel subset $\gM$ of $(0,+\infty)$. Then there is a 
$\U_q(gl_2)$-invariant linear functional $h_{\mu_0}$ on the 
$\ast$-subalgebra $\cO (SU_q(2))\otimes C_c(0,+\infty)$ of 
$\hat{\cO}_e(\C^2_q)$ such that
$$
h_{\mu_0}(x\varphi(\cR_q))=h(x){\int^\infty_0}\varphi(t) 
d\mu (t),~x\in\cO(SU_q(2)),\varphi\in C_c(0,+\infty).
$$
For instance, if $\mu_0$ is a Dirac measure 
$\delta_{t_0}$ , then $\mu$ is supported on the 
points $t_0p^n$, $n\in\Z$, and $h_{\mu_0}(\varphi(\cR_q))$ is given by the 
Jackson integral
$$
h_{\mu_0}(\varphi(\cR_q))={\int^\infty_0}\varphi(t)d\mu(t)=
{\sum^{+\infty}_{n=0}}\varphi(t_0p^n)p^n.
$$

\subsection{Heisenberg representations of cross product algebras}
Let $\U_0\lti\X$ be a right cross product algebra as in Section 2. It is 
well known that there exists a unique homomorphism $\pi$ of 
$\U_0\lti\X$ into the algebra $L(\X)$ of linear mappings of $\X$ such that $\pi(x)y=xy$ and $\pi(f)y=y\anf S^{-1}(f)$ for $x,y\in\X$ and $f\in\U_0$. 
(Indeed, one easily checks that $\pi(x)\pi(f)=
\pi(f_{(1)})\pi(x\anf f_{(2)})$, so there is a well defined (!) homomorphism $\pi$ with these 
properties.) In this subsection we develop a Hilbert space 
version of this algebraic fact.

Suppose $\X$ is a $\ast$-algebra with unit element. Let $h$ be a state 
on $\X$, i.e. $h$ is a linear functional on $\X$ such that 
$h(x^\ast x)\ge 0$ for $x\in\X$ and $h(1)=1$, and let $\pi_h$ denote the 
$GNS$-representation $\pi_h$ of $\X$ (see e.g.\ [S], Section 8.6). By the definition 
of the $GNS$-representation, there is a vector $v_h$ in the domain of $\pi_h$ 
such that $\D_h=\pi_h(\X)v_h$ is dense in the underlying Hilbert space 
$\Hh_h$, $\pi_h(x) (\pi_h(y)v_h)=\pi_h (xy) v_h$ for $x,y\in\X$ and 
\begin{equation}\label{gnsx}
h(x)=\langle \pi_h(x)v_h,v_h\rangle,~ x\in\X.
\end{equation}

\bn
{\bf Proposition 5.3}
{\it Let $\X$ be a right $\U$-module $\ast$-algebra of a Hopf 
$\ast$-algebra $\U$ and let $\U_0$ be a unital $\ast$-subalgebra and a 
right coideal of $\U$. Suppose that $h$ is a $\U_0$-invariant state on $\X$. Then there exists a unique $\ast$-representation $\tilde{\pi}_h$ 
of the $\ast$-algebra $\U_0\lti\X$ on the domain $\D_h=\pi_h(\X)v_h$ such 
that $\tilde{\pi}_h\lceil\X$ is the $GNS$-representation of $\X$ and 
$\tilde{\pi}_h(f) v_h=\varepsilon (f) v_h$ for $f\in\U_0$. For $f\in\U_0$ 
and $x\in\X$, we have}
\begin{equation}\label{gnsf}
\tilde{\pi}_h(f)(\pi_h(x)v_h)=\pi_h(x\anf S^{-1} (f))v_h.
\end{equation}

\mn
{\bf Proof.} First we show that the linear mapping $\tilde{\pi}_h(f)$ 
given by (\ref{gnsf}) is well defined. Let $x\in\X$ be such that 
$\pi_h(x)v_h=0$ and let $f,g\in\U_0$. Using (\ref{gnsx}), (\ref{modstar}), and the 
$\U_0$-invariance of $h$ we obtain
\begin{align*}
&\langle \pi_h(x\anf S^{-1}(f))v_h, \pi_h(x\anf S^{-1} (g))v_h\rangle
=h((x\anf S^{-1}(g))^\ast(x\anf S^{-1}(f)))\\
&=h((x^\ast\anf g^\ast)(x\anf 
S^{-1}(f_{(2)})))\varepsilon(f_{(1)})
=h((x^\ast\anf g^\ast)(x\anf S^{-1}(f_{(2)}))\anf f_{(1)})\\
&=h((x^\ast\anf g^\ast f_{(1)})(x\anf S^{-1}(f_{(3)}) f_{(2)}))
=h((x^\ast\anf g^\ast f)x)\\
&=\langle \pi_h(x)v_h,\pi_h((x^\ast\anf 
g^\ast f)^\ast)v_h\rangle =0.
\end{align*}
Setting $f=g$, we get $\pi_h(x\anf S^{-1}(f))v_h=0$. Hence the linear 
mapping $\tilde{\pi}_h(f)$ is well defined by (\ref{gnsf}).

Obviously, the map $f\rightarrow\tilde{\pi}_h(f)$ is an algebra 
homomorphism of $\U_0$ into $L(\D_h)$. From (\ref{gnsf}) and the 
definition of the $GNS$-representation $\pi_h$ of $\X$ it follows 
easily that $\tilde{\pi}_h$ extends to an algebra homomorphism, 
denoted again by  $\tilde{\pi}_h$, of $\U_0\lti\X$ into $L(\D_h)$ 
such that $\tilde{\pi}_h\lceil \X=\pi_h$. Next we show that 
$\tilde{\pi}_h$ preserves the involution. Let $f\in \U_0$ and 
$x,y\in\X$. Using (\ref{gnsx}), (\ref{gnsf}), (\ref{modstar}), and 
the $\U_0$-invariance of $h$ we compute
\begin{align*}
&\langle \tilde{\pi}_h(f)\pi_h(x) v_h,\pi_h(y)v_h\rangle=
\langle\pi_h(y^\ast(x\anf S^{-1}(f)))v_h,v_h\rangle
=h(y^\ast(x\anf S^{-1}(f)))\\
&=h(y^\ast(x\anf S^{-1}(f_{(2)})))\varepsilon(f_{(1)})
=h((y^\ast(x\anf S^{-1}(f_{(2)}))\anf f_{(1)})\\
&=h((y^\ast\anf f_{(1)})(x\anf S^{-1} (f_{(3)}) f_{(2)}))
=h((y^\ast\anf f)x)=h((y\anf S^{-1} (f^\ast))^\ast x)\\
&=\langle \pi_h((y\anf S^{-1} (f^\ast))^\ast x)v_h,v_h\rangle
=\langle \pi_h(x) v_h,\tilde{\pi}_h(f^\ast)\pi_h(y) v_h\rangle.
\end{align*}
This shows that $\tilde{\pi}_h\lceil \U_0$ is a $\ast$-representation of 
the $\ast$-algebra $\U_0$. The restriction of $\tilde{\pi}_h$ to $\X$ is 
the GNS-representation and so a $\ast$-representation. Since the 
$\ast$-algebra $\U_0\lti\X$ is generated by $\U_0$ and $\X$ (by Lemma 2.1), 
$\tilde{\pi}_h$ is a $\ast$-representation of the $\ast$-algebra 
$\U_0\lti\X$ on the domain $\D_h$.

Finally, we prove the uniqueness assertion. Let $\pi$ be a representation of 
$\U_0\lti\X$ such that $\pi\lceil\X=\pi_h$ and $\pi(f) v_h=
\varepsilon(f) v_h$. Using (\ref{crossfx}) and the fact that $\pi$ is an 
algebra homomorphism we obtain
\begin{align*}
&\pi (f) \pi_h(x)v_h=\pi (fx)v_h=\pi((x\anf S^{-1} (f_{(2)})f_{(1)})v_h\\
&=\pi(x\anf S^{-1}(f_{(2)}))\pi(f_{(1)})v_h=\pi(x\anf S^{-1}(f))v_h
=\tilde{\pi}_h(f)\pi_h(x) v_h
\end{align*}
for $f\in\U_0$ and $x\in\X$. That is, we have $\pi=\tilde{\pi}_h$.\hfill 
$\Box$

\sn
Let us denote the closure of the 
$\ast$-representation $\tilde{\pi}_h$ of $\U_0\lti \X$ 
by $\pi_h$. 

\mn
{\bf Definition 5.4} The closed $\ast$-representation $\pi_h$ of $\U_0\lti\X$ 
defined above is called the {\it Heisenberg representation} 
or the {\it left regular representation} of the cross 
product $\ast$-algebra $\U_0\lti \X$ with respect to the $\U_0$-invariant 
state $h$.

\subsection{A uniqueness theorem for the Heisenberg 
representation of ${\bf \U_q(su_2)\lti \cO(SU_q(2))}$ }
By Lemma 5.2 the $\U_q(su_2)$- module $\ast$-algebra $\cO(SU_q(2))$ has a unique $\U_q(su_2)$-invariant state. Hence the cross product algebra ${\U_q(su_2){\lti}\cO(SU_q(2))}$  has a unique Heisenberg representation $\pi_h$.

\mn
{\bf Theorem 5.5} (i) {\it 
Let $\pi$ be a closed $\ast$-representation of 
$\U_q(su_2)\lti\cO(SU_q(2))$ such that the restriction of $\pi$ to $\U_q(su_2)$ is integrable.
Then $\pi$ is unitarily equivalent to a direct sum of Heisenberg 
representations $\pi_h$.}\\
(ii) {\it The Heisenberg representation $\pi_h$ of $\U_q(su_2)\lti\cO(SU_q(2))$ is 
irreducible and its restriction to $\U_q(su_2)$ is an integrable 
representation.}

\mn
An immediate consequence of this theorem is

\mn
{\bf Corollary 5.6} {\it Let $\pi$ be a closed $\ast$-representation of 
$\U_q(su_2)\lti \cO(SU_q(2))$. Then $\pi$ is unitarily equivalent to the Heisenberg 
representation $\pi_h$ if and only if $\pi$ is 
irreducible and the restriction of $\pi$ to $\U_q(su_2)$ is integrable.}

\mn
{\bf Proof of Theorem 5.5.} 
(i): The crucial step of this proof is to show that there exists a 
non-zero vector $v_0$ in the domain $\D$ of the representation $\pi$ such that
\begin{equation}\label{weight}
\pi(E)v_0=\pi(F)v_0=0~{\rm and}~\pi(K)v_0=v_0.
\end{equation}
By assumption,  $\pi \lceil \U_q(su_2)$ is a direct sum of representations 
$T_{l_i}, i\in I$, 
where $l_i\in {\frac{1}{2}}\N_0$. Let us take a non-zero lowest weight 
vector $v_l\in\D$ of  weight $-l$ in the representation space of one of 
these representations $T_l,l=l_i,i\in I$. That is, we have $\pi(F)v_l=0$ and 
$\pi(K)v_l=q^{-l} v_l$. If $l=0$, then $T_l$ is the trivial 
subrepresentation of $\U_q(su_2)$ and we are done. If $l\ge 1/2$, then we put
\begin{align*}
&w_l := \pi(c) \pi(E) v_l + q^{-l-3/2} [2l]_q \pi(a) v_l,\\
&w^\prime_l := \pi(d) \pi(E) v_l + q^{-l-3/2} [2l]_q \pi(b) v_l.
\end{align*}
Using the cross relations of the algebra  $\U_q(su_2)\lti\cO(SU_q(2))$ we 
derive  
\begin{equation*}
\pi(F) w_l = \pi (F) w^\prime_l = 0,\pi(K) w_l = q^{-(l-1/2)} w_l, 
\pi (K) w^\prime_l = q^{-(l-1/2)} w^\prime_l.
\end{equation*}
Further, we have $\pi(a)w^\prime_l-\pi(b)w_l=\pi(ad-qbc)\pi(E)v_l=
\pi(E)v_l\ne 0$, because we assumed that $l\ge 1/2$. Hence at least one 
of the vectors $w_l,w_{l^\prime}$ is non-zero.
Thus we have constructed a non-zero lowest weight vector of weight ${-}(l{-}1/2)$. Proceeding by induction we obtain after $2l$ steps a non-zero lowest weight vector $v_0\in\D$ of weight $0$, 
that is, $\pi(F)v_0=0$ and $\pi(K) v_0=v_0$. From the cross 
relations of $\U_q(su_2)\lti \cO(SU_q(2))$ we compute $\pi(E)w_{1/2}=
\pi(E)w^\prime_{1/2}=0$, so the vector $v_0$ satisfies also $\pi(E)v_0=0$.

After norming we can assume that ${\parallel v_0\parallel}=1$. Obviously, 
(\ref{weight}) implies 
\begin{equation}\label{invfv}
\pi(f)v_0=\varepsilon (f) v_0, f\in\U_q(su_2).
\end{equation}
Define $h(x)=\langle \pi (x)v_{0},v_0\rangle$ for $x\in\cO (SU_q(2))$. 
Using (\ref{adright}) and (\ref{gnsx}) we obtain
\begin{align*}
&h(x\anf f)=\langle\pi(ad_R(f)x)v_0,v_0\rangle=\langle 
\pi(S(f_{(1)}))\pi(x)\pi(f_{(2)})v_0,v_0\rangle\\
&=\langle\pi(x)\pi(f_{(2)})v_0,\pi(S(f_{(1)})^\ast)v_0\rangle=\\
&=\langle \pi(x)\varepsilon(f_{(2)})v_0, 
\varepsilon(S(f_{(1)})^\ast)v_0\rangle =\varepsilon (f)
\langle \pi(x)v_0,v_0\rangle =\varepsilon (f)h(x).
\end{align*}
That is, $h$ is a $\U_q(su_2)$-invariant linear functional on 
$\cO (SU_q(2))$ such that $h(1)=1$. By 
Lemma 5.2,
$h$ is the Haar state of $\cO(SU_q(2))$.
Put $\D_0:=\pi(\cO(SU_q(2))v_0$. Let $f\in\U_q(su_2)$ and 
$x\in\cO(SU_q(2))$. By (\ref{crossfx}) and (\ref{invfv}), we have
\begin{align}\label{invfd}
&\pi(f)\pi(x) v_0=\pi(fx)v_0=\pi((x\anf S^{-1}(f_{(2)}))\pi(f_{(1)}) 
v_0\nonumber\\
&=\pi(x\anf S^{-1}(f_{(2)}))\varepsilon (f_{(1)})v_0=\pi(x\anf S^{-1}(f))v_0,
\end{align}
so $\pi(f)$ leaves the domain $\D_0$ invariant. Hence $\D_0$ is invariant 
under the representation $\pi$ of the whole algebra $\U_q(su_2)\lti\cO(SU_q(2))$. 
By (\ref{invfv}), we have $\pi(f)v_0=\varepsilon(f)v_0$. 
Since $h(x)=\langle \pi(x) v_0,v_0\rangle$ for 
$x\in\cO(SU_q(2))$, the restriction of $\pi\lceil 
\cO(SU_q(2))$ to the domain $\D_0$ is unitarily equivalent to the  
$GNS$-representation $\pi_h$ on $\D_h=\pi_h(\cO(SU_q(2))v_h$. Let 
$\pi_0$ be the closure of the restriction of the $\ast$-representation 
$\pi$ of $\U_q(su_2)\lti\cO(SU_q(2))$ to $\D_0$. By the preceding and the 
uniqueness assertion of Proposition 5.3, $\pi_0$ is unitarily equivalent to 
the Heisenberg representation $\pi_h$.

By assumption the representation $\pi\lceil\U_q(su_2)$ is a 
direct sum of finite dimensional representations $T_{l_i}, 
l_i\in{\frac{1}{2}}\N_0$. Further, the operators $\pi(x)$, $x\in\cO(SU_q(2))$, 
are bounded. From these facts it follows that $\pi$ decomposes into a direct 
sum $\pi_0\oplus\pi^\prime$, where $\pi^\prime$ satisfies again the 
assumptions of (i). A standard argument based on Zorn's Lemma gives the 
assertion.

(ii): Let $t^{(l)}_{ij}$ be the matrix elements of the 
spin $l$ corepresentation of $SU_q(2)$ and $V_{lj}:= \Lin\{t^{(l)}_{ij}; 
i=-l,{\dots},l\}$. By (\ref{rightact}) and (\ref{gnsf}), 
$V_{lj}$ is invariant under $\pi_h(\U_q(su_2))$, 
and the subrepresentation of 
$\U_q(su_2)$ on $V_{lj}$ is unitarily 
equivalent to $T_l$, $l\in{\frac{1}{2}}\N_0$. Since 
the set $\{t^{(l)}_{ij};i,j=-l,{\dots},l,l\in{\frac{1}{2}}\N_0\}$ is a
 vector space basis of $\cO(SU_q(2))$ by the Peter-Weyl theorem, 
$\pi_h\lceil\U_q(su_2)$ is a direct 
sum of representations $T_l, l\in{\frac{1}{2}}\N_0$. That is,
$\pi_h\lceil \U_q(su_2)$ is integrable. 

Finally, we prove that $\pi_h$ is irreducible. Assume to the contrary that 
$\pi_h$ is the direct sum of non-trivial representations $\pi_1$ and 
$\pi_2$.
Since $\pi_h\lceil \U_q(su_2)$ 
is integrable as just shown, 
$\pi_j\lceil \U_q(su_2)$, $j{=}1,2$, is also integrable. By (i), $\pi_j$ is a direct sum of
Heisenberg representations. Hence $\pi_h$ is unitarily equivalent to 
a sum of more than two copies of $\pi_h$. 
Then the dimension of the space of $\U_q(su_2)$-invariant vectors 
in $\pi_h(\cO(SU_q(2))v_h$ is larger than one which contradicts the 
uniqueness of the Haar functional of $SU_q(2)$.\hfill $\Box$

\subsection{A quantum trace formula for the Haar state of 
${\bf \cO(SU_q(2))}$}
If $T$ is a finite dimensional representation of $\U_q(sl_2)$, then 
the {\it quantum trace} ${\rm Tr}_q(f):={\rm Tr}~ T(K^{-2}f)$, 
$f\in\U_q(sl_2)$, is an $ad_R$-invariant linear functional on $\U_q(sl_2)$, that is, ${\rm Tr}_q (\ad_R(f)g)=\varepsilon (f) {\rm Tr}_q (g)$ for $f,g\in\U_q(sl_2)$. This well known fact is based on the trace 
property and the relation $K^{-2} S(f)=S^{-1}(f) K^{-2}$. For representations on 
infinite dimensional Hilbert spaces the quantum trace in the above form 
does not make sense. In this subsection and in Subsection 6.3 
below we develop two variants of the quantum trace and describe the Haar state of 
$SU_q(2)$ in this manner.

For $z\in\C$, ${\rm Re}~ z>1$, we define a holomorphic function
$$
\zeta(z)= (q^{-1}-q)^{2z-1} 
{\sum\limits^\infty_{n=1}} n 
(q^{-\frac{n}{2}}-q^{\frac{n}{2}})^{-2z}(q^{-n} - q^n).
$$
Note that $q^{-1}-q>0$ and $q^{-\frac{n}{2}}-q^{\frac{n}{2}}>0$ for $n\in\N$, 
since $0<q<1$. Let $\cC$ denote the closure of the operator 
$\pi_h(C_q)$, where $C_q$ is the Casimir element (\ref{casi}) of 
$\U_q(su_2)$. Recall that $h$ is the Haar state of $\cO(SU_q(2))$.

\mn
{\bf Theorem 5.7} {\it If $z\in\C$, $Re~z>1$, and $x\in\cO(SU_q(2))$, then the 
closure of the operator $\cC^{-z}\pi_h(K^{-2}x)$ is of trace class and 
we have}
\begin{equation}\label{hcas}
h(x)=\zeta(z)^{-1}~ {\rm Tr}~
 {\overline{\cC^{-z}\pi_h(K^{-2}x)}}~.
\end{equation} 

\mn
{\bf Proof.} 
The representation $\pi_h\lceil\U_q(su_2)$ is the 
direct sum of representations $(2l+1)T_l$, $l\in \frac{1}{2}\N_0$. Hence, by (\ref{tlkef}) and (\ref{tlc}), the operators  ${\overline{\pi_h(K)}}$, $|{\overline{\pi_h(E)}}|$, $|{\overline{\pi_h(F)}}|$, and 
$\cC={\overline{\pi_h(C_q)}}$ have a common orthonormal basis of 
eigenvectors with eigenvalues $q^j$, $\alpha_{j+1,l}$, 
$\alpha_{jl}$, $j=-l,{\dots},l$, each with multiplicity $2l+1$, and $[l+{1/2}]^2_q$ with multiplicity 
$(2l{+}1)^2$, respectively. Since 
 $|\alpha_{jl}|\le {\rm const.}~q^{-l}$, $[l+1/2]^{-2}_q\le
{\rm const.}~q^{2l}$ and $\pi_h(x)$, $x\in\cO(SU_q(2))$, is bounded, 
it follows  that the closures of  the operators $\cC^{-z}\pi_h(K^{-2} x)$, 
$\cC^{-z}\pi_h(K^{-1} Ex)$ and $\cC^{-z}\pi_h(K^{-1} Fx)$ are of trace 
 class when ${\rm Re}~ z>1.$ Thus, $h_z(x):={\rm Tr}~ 
{\overline{\cC^{-z} \pi_h(K^{-2}x)}}$ is well defined for 
$x\in\cO(SU_q(2))$ and $z\in\C$, ${\rm  Re}~ z>1$.
By the cross relations of $\U_q (su_2)\lti\cO(SU_q(2))$, each element 
$xK^{-1} E$ is of 
the form $K^{-1} Ex^\prime +K^{-2}x^{\prime\prime}$ with 
$x^\prime,x^{\prime\prime}\in\cO(SU_q(2))$. Hence the closures of 
$\cC^{-z}\pi_h(xK^{-1}E)$ and likewise of 
$\cC^{-z}\pi_h(xK^{-1}F)$ are also of trace class.

Suppose $z\in\C$, ${\rm Re}~ z>1$, and $x\in\cO(SU_q(2))$. There is 
$y\in\cO(SU_q(2))$ such that $xK=Ky$. Using the facts of the 
preceding paragraph we conclude 
\begin{align*}
&{\rm Tr}~ {\overline{\cC^{-2z}\pi_h(yK^{-1} E)}}=
{\rm Tr}~\cC^{-z}{\overline{\cC^{-z}\pi_h(yK^{-1} E)}}\\
&={\rm Tr}~ {\overline{\cC^{-z}\pi_h(yK^{-1} E)}}\cC^{-z}=
{\rm Tr}~\cC^{-z}{\overline{\pi_h(y)}}{\overline{\pi_h(K^{-1} E)
\cC^{-z}}}\\
&={\rm Tr}~\pi_h({\overline{K^{-1} E)\cC^{-z}}}
\cC^{-z}\pi_h(y)=
{\rm Tr}~{\overline{\cC^{-2z}\pi_h(K^{-1} Ey)}}.\\
\text{From (\ref{adright}) and}~&\text{the latter formula we get}\hspace{3cm}\\
&h_{2z}(x\anf E)=h_{2z}(ad_R(E)x)={\rm Tr}~{\overline{\cC^{-2z}
\pi_h(K^{-2}(KxE-qExK))}}\\
&={\rm Tr}~ {\overline{\cC^{-2z}\pi_h(yK^{-1}E-K^{-1}Ey)}}\\
&={\rm Tr}~{\overline{\cC^{-2z}\pi_h(yK^{-1}E)}}-
{\rm Tr}~{\overline{\cC^{-2z}\pi_h(K^{-1}Ey)}}=0.
\end{align*}
Similarly, one shows that $h_{2z}(x\anf F)=0$ and $h_{2z}(x\anf K^{\pm 1})=h_{2z}(x)$. 
Hence $h_{2z}(x\anf f)=\varepsilon(f)h_{2z} (x)$ for $f{\in}\U_q(su_2)$. Thus $h_{2z}$ is a $\U_q(su_2)$-invariant linear functional on 
$\cO(SU_q(2))$. By Lemma 5.2, $h_{2z}(1)^{-1} h_{2z}$ is the Haar state $h$ 
of $\cO(SU_q(2))$.

Now we compute $h_z(1)$ for $z\in\C$, ${\rm  Re}~ z>1$. The trace of the 
restriction of the operator $\cC^{-z}\pi_h(K^{-2})$ to the invariant 
subspace $V_l$ is 
$$
{\sum\limits^l_{i,j=-l}} [l+{1/2}]^{-2z}_q q^{-2j}= (2l+1)[l+
{1/2}]^{-2z}_q [2l+1]_q.
$$
Summing over $l\in\frac{1}{2}\N_0$ or equivalently over $n=2l{+}1\in\N$, we get 
$$h_z(1)={\rm Tr}~{\overline{\cC^{-z}\pi_h(K^{-2})}}=
{\sum\limits^\infty_{n=1}} n [{n/2}]^{-2z}_q [n]_q=\zeta(z).
$$
By the preceding we have proved that 
$h=h_z(1)^{-1}h_z=\zeta(z)^{-1}h_z$ for 
$z\in\C$, ${\rm Re}~ z>2$. Let $x \in \cO(SU_q(2))$. Then, 
$\zeta(z)h(x)$ and $h_z(x)$ are holomorphic functions of 
$z\in\C$, ${\rm Re}~ z>1$. As just shown, the two functions are equal for 
${\rm Re}~ z>2$. Hence they coincide also for ${\rm Re}~ z>1$. \hfill $\Box$\\
 
The function $\zeta(z)$ is called the {\it zeta function} of the
quantum group $SU_q(2)$.
The functional $h_z$ and so formula (\ref{hcas}) can be rewritten as
$$
h_z(x)={\rm Tr}~{\overline{\pi_h(x)}}~
 {\overline{\pi_h(K^{-2})\pi_h(C_q)^{-z}}}=
{\rm Tr}~ {\overline{\pi_h(C_q)^{-z}\pi_h(K^{-2})}}~ {\overline{\pi_h(x)}}.
$$
Note that the operators ${\overline{\pi_h(K^{-2})\pi_h(C_q)^{-z}}}$ and 
${\overline{\pi_h(C_q)^{-z}\pi_h(K^{-2})}}$ are of trace class 
if ${\rm Re}~ z>1$ and ${\overline{\pi_h(x)}},x\in\cO(SU_q(2))$, is bounded.

\section{Representations of cross product $\ast$-algebras}  

We now develop the second approach to representations of cross 
product algebras. That is, we begin with a $\ast$-representation 
of one of the coordinate $\ast$-algebras $\cO(SU_q(2))$, ${\hat\cO}(\C^2_q)$, 
and $\cO(\R^3_q)$ as described in 4.2 and try to 
complete it to a $\ast$-representation 
of the cross product $\ast$-algebra. In doing so we require that 
certain relations derived by formal algebraic manipulations hold in the 
operator-theoretic sense in the Hilbert space. For notational simplicity 
we suppress the representation and write $x$ for $\pi(x)$ when no 
confusion can arise. Recall that $0<q<1$, $p>0$ and $p\ne 1$.
 
\mn
\subsection{Representations of the ${\bf \ast}$-algebra ${\bf \U_0\lti\cO(SU_q(2))}$}
Suppose we have a $\ast$-representation of $\U_0\lti\cO(SU_q(2))$ on a Hilbert 
space such that its restriction to $\cO(SU_q(2))$ 
is of the form described in 4.2.

We assume that there exist dense linear subspaces $\E$ and $\D_0$ of $\G$ 
and $\Hh_0$, respectively, such that $v\E=\E, w\D_0=\D_0$ and 
$\E\oplus\D$ is invariant under the $X_j$, $j=0,1,2$, where 
$\D=\Lin\{\eta_n; \eta \in \D_0,~n\in\N_0\}$.

\sn
{\bf 1. Step:} 
First we show that $\G=\{0\}$. Since $X^\ast_2=X_0$ and $b=c=0$ on $\G$, 
it follows from the relations $bX_2=q^{-1}X_2b$ and $c X_0{=}qX_0c$ that 
$X_2$ leaves the subspace $\G$ invariant. Thus, 
$cX_2\varphi{=}0$ for $\varphi\in\E$. Therefore, 
since $c X_2{=}q X_2 c{+} a$, we obtain $a \varphi {=} 0$ for 
$\varphi \in \E$. Because $a {=} v$ is unitary on $\G$ and 
$\E$ is dense in $\G$, the latter implies that $\G{=}\{0\}$. 

Since $\Hh_0=\ker a$, it follows from the relation 
$a Y_1=q^{-2}Y_1a$ that $Y_1$ leaves $\Hh_0$ invariant, 
so there is a symmetric operator $G_0$ on $\Hh_0$ such 
that $Y_1\eta_0{=}G_0\eta_0,\eta_0\in\D_0$. From 
$d^n Y_1{=}q^{2n}Y_1d^n$ we conclude that $Y_1\eta_n{=}q^{-2n} G_0\eta_n$. 
The relation $c Y_1{=}q^2Y_1 c$ implies that 
$w^\ast G_0\eta_n{=}q^2G_0w^\ast \eta_n$.

From the relations $aX_2=q^{-1}X_2a$ and $d X_2=qX_2 d+b$ it follows by induction on $n$ that 
$X_2$ maps $\Hh_n$ into $\Hh_n+\Hh_{n-1}$. Hence there exist linear 
operators $T_n$ and $R_n$ on the Hilbert space $\Hh_0$ such that 
$X_2\eta_n=T_n\eta_n+R_n\eta_{n-1}, \eta\in\D_0$. Inserting this into the relation $a X_2=q^{-1}X_2a$ shows that $T_n=q^{-1} T_{n-1}$. Thus we obtain $T_n=q^{-n} T_0$. Comparing the $(n{-}1)$th  components of the relation $bc X_2=X_2bc+ba$ gives $R_n=q^{-n}\lambda_n\lambda^{-1}w$. Since $b X_2=q^{-1}X_2b$, we obtain $w T_0=q^{-1} T_0 w$. 

We have not yet used the commutation relations (\ref{tan1})--(\ref{tan3}) 
of the generators $X_0,X_2,X_1$. Equation (\ref{tan2}) is equivalent to 
$Y_1X_2=q^4X_2Y_1$. Using the formulas for $Y_1$ and $X_2$ obtained 
above we get $G_0T_0=q^4T_0G_0$. Since $X_0 =X^\ast_2$, we 
have $X_0\eta_n=T^\ast_n\eta_n+R^\ast_{n+1}\eta_{n+1},\eta\in\D_0$. 
Inserting the above expressions of $X_0,X_2,T_n,R_n$ 
into relation (\ref{tan3}) yields the equation $qT_0T^\ast_0-q^{-1} T^\ast_0T_0=\lambda^{-1} G_0$. 

Summarizing the first step, we have shown that  
$X_0,X_2, Y_1$ act as
\begin{align}
\label{x2}
&X_2\eta_n=q^{-n} T_0\eta_n+q^{-n}\lambda_n\lambda^{-1}w \eta_{n-1},\\
\label{x0}
&X_0\eta_n=q^{-n} T^\ast_0\eta_n+
q^{-n-1}\lambda_{n+1}\lambda^{-1} w^\ast\eta_{n+1},\\
\label{x1}
&Y_1\eta_n=q^{-2n} G_0\eta_n,\eta\in\D_0,
\end{align}
where the operators $T_0,G_0,w$  satisfy the consistency conditions
\begin{align}
\label{cons1}
&w G_0w^\ast=q^{-2}G_0,~G_0T_0=q^4T_0G_0,~ wT_0w^\ast=q^{-1}T_0,\\
\label{cons2}
&qT_0T^\ast_0-q^{-1} T^\ast_0T_0=\lambda^{-1} G_0
\end{align}
on the domain $\D_0$ of $\Hh_0$. Conversely, if operators $T_0,G_0,w$ on 
the Hilbert space $\Hh_0$ are given such that $G_0$ is symmetric, 
$w$ is unitary and relations (\ref{cons1}) and (\ref{cons2}) are 
fulfilled on a dense domain $\D_0$ of $\Hh_0$ which is invariant 
for the operators $T_0, T^\ast_0,G_0,w,$ and $w^\ast$, then the
formulas (\ref{opsu2}), (\ref{x2})--(\ref{x1}) define a 
$\ast$-representation of $\U_0\lti\cO(SU_q(2))$. Indeed, by 
straightforward computations it can be checked that the 
defining relations of $\U_0\lti\cO(SU_q(2))$ are satisfied.

\sn
{\bf 2. Step:}
Next we analyze the triple of operators $T_0,G_0,w$ on the Hilbert 
space $\Hh_0$ satisfying the consistency conditions 
(\ref{cons1})--(\ref{cons2}). 

We suppose that the closure of the operator $Y_1$ is self-adjoint. 
By the above assumption on the domain, we can suppose that 
$G_0$ is self-adjoint. Let $g(\lambda),\lambda\in\R$, denote the 
spectral projections of $G_0$. The Hilbert space $\Hh_0$ decomposes 
into a direct sum 
\begin{equation}
\label{dir}
\Hh_0=g((-\infty,0))\Hh_0\oplus g(\{0\})\Hh_0\oplus g((0,\infty))\Hh_0
\end{equation}
of reducing subspaces of $G_0$ where $G_0{<} 0$, $G_0=0$ and $G_0{>}0$, 
respectively. Since $wG_0 w^\ast=q^{-2} G_0$, the direct sum (\ref{dir}) 
reduces $w$ as well. Since $G_0 T_0 = q^4 T_0 G_0$, we assume that (\ref{dir}) reduces also $T_0$. 
Thus we are lead to study relations (\ref{cons1})--(\ref{cons2}) in the
three cases $G_0=0$, $G_0 {<} 0$, $G_0{>} 0$ separately.

\sn
{\it Case I.} $G_0=0$

Then relations (\ref{cons1})--(\ref{cons2}) read $wT_0w^\ast=q^{-1}T_0$ and $q^2T_0T^\ast_0=T^\ast_0T_0$. Obviously, 
there is the trivial representation where $T_0=0$ and $w$ 
is arbitrary. 
Since $\ker T_0=\ker T^\ast_0$ is reducing for 
$T_0$ and $w$, we can consider the cases $T_0=0$ and $\ker T_0=\{0\}$ separately.
Assume now that $\ker T_0=\{0\}$ and let 
$T_0=v_0|T_0|$ be the polar decomposition of $T_0$. 
Since $\ker T_0=\ker T^\ast_0=\{0\}, v_0$ is unitary. From  
$wT_0w^\ast=q^{-1} T_0$ we get $w |T_0| w^\ast = q^{-1} |T_0|$. The 
relation $q^2 T_0T^\ast_0=T^\ast_0 T_0$ implies that 
$q^2 v_0 |T_0|^2 v^\ast_0=|T_0|^2$ and so 
$q v_0 |T_0| v^\ast_0= |T_0|$. Hence we have  
$w^\ast v_0 |T_0|v_0^\ast w=q^{-1} w^\ast |T_0| w=|T_0| $ . 
Using the preceding relations and the fact that $w$ and $v_0$ are unitary, 
we get $w^\ast T_0=w^\ast v_0|T_0|=|T_0|w^\ast v_0=
q T_0 w^\ast= qv_0|T_0|w^\ast=v_0 w^\ast|T_0|$. Since 
$\ker |T_0|=\ker T_0=\{0\}$, we obtain $w^\ast v_0=v_0w^\ast$ and so 
$wv_0=v_0w$.

Since $w$ is unitary, by Lemma 4.2(i) the relation $w|T_0|w^\ast=q^{-1}|T_0|$ leads to a 
representation $w\zeta_k=\zeta_{k+1}$, $|T_0|\zeta_k=q^k A_{00}\zeta_k$ 
on $\Hh_0=\displaystyle\mathop{\oplus}\nolimits^{\infty}_{n=-\infty} \Hh_{0k}$, 
where $\Hh_{0k}=\Hh_{00}$ and $A_{00}$ is a self-adjoint operator on the 
Hilbert space $\Hh_{00}$ such that $\sigma(A_{00})\sqsubseteq (q,1]$. 
The operator $w^\ast v_0$ commutes with $|T_0|$ and so with the spectral 
projections of $|T_0|$. Since $\sigma (A_{00})\sqsubseteq (q,1]$, 
$w^\ast v_0$ leaves each space $\Hh_{0k}$ invariant. Therefore, 
since $wv_0=v_0w$, there is a unitary operator $v_{00}$ on $\Hh_{0k}$ such that $v_0\zeta_k=v_{00}\zeta_{k+1}, k\in\Z$. Using
again the relation $q^2T_0T^\ast_0=T^\ast_0T_0$, it follows that the 
operator $N:=v_{00} A_{00}$ on $\Hh_{00}$ is normal. Thus we have 
$T_0\zeta_k=q^k N\zeta_{k+1}$ with  $N$ normal. This completes the 
treatment of Case I. 

\mn
Next we treat the cases $G_0<0$ and $G_0>0$. Let us set $G_0=\delta H^2_0$ and 
$H_0:=\epsilon |G_0|^{1/2}$, where $\epsilon,\delta\in\{1,-1\}$. Since 
$G_0T_0=q^4T_0G_0$, it is natural to assume that  $H_0T_0=q^2T_0H_0$. 
(It can be shown that the relation $H_0T_0=-q^2 T_0H_0$ does not have a non-trivial 
solution for $H_0>0$.) Set $S_0:=H_0^{-1}T_0$. We rewrite the consistency 
conditions in terms of $H_0, S_0, w$ by formal (!) algebraic 
manipulations. Since $S_0=H^{-1}_0T_0=q^{-2}T_0H^{-1}_0$ and so 
$S^\ast_0=T^\ast_0 H^{-1}_0=q^{-2} H^{-1}_0 T^\ast_0$, (\ref{cons2}) is 
formally equivalent to
\begin{equation}\label{cons22}
S_0 S^\ast_0-q^2 S^\ast_0 S_0=-\delta(1{-}q^2)^{-1}~,
\end{equation}
where $\delta\in\{1,-1\}$. The three relations (\ref{cons1}) can be rewritten as 
\begin{equation}\label{cons21}
w H_0w^\ast=q^{-1} H_0,~ H_0S_0=q^2 S_0 H_0,~ w S_0w^\ast=S_0~.
\end{equation}
We now solve the relations (\ref{cons22}) and (\ref{cons21}) in a rigorous 
manner.

\mn
{\it Case II.} $G_0 < 0~ (\delta=-1)$

Since $0 <q<1$ and $-\delta(1{-}q^2)^{-1}>0$, relation (\ref{cons22}) has 
three series of representations by Lemma 4.3. Let us begin with the {\it Fock representation}. 
Then we have 
\begin{equation}\label{so1}
S_0\zeta_k=(1{-}q^2)^{-1}\lambda_k\zeta_{k-1},~ S^\ast_0\zeta_k=(1{-}q^2)^{-1}\lambda_{k+1}\zeta_{k+1}
\end{equation}
acting on the direct sum Hilbert space 
$\Hh_0{=} \displaystyle\mathop{\oplus}\nolimits^\infty_{k=0} \Hh_{0k}$, 
where $\Hh_{0k}{=}\Hh_{00}$. Since 
$w S_0w^\ast{=} S_0$ by (\ref{cons21}), $w$ commutes with the spectral 
projections of $S_0S^\ast_0$ and so $w$ leaves each space $\Hh_{0k}$ 
invariant. Hence, by (\ref{so1}), the relation $w S_0w^\ast=S_0$ implies that 
there is a unitary $w_0$ on $\Hh_{00}$ such that $w\zeta_k=w_0\zeta_k$. 
The relation $H_0 S_0=q^2 S_0H_0$ yields $H_0S^\ast_0S_0=S^\ast_0S_0H_0$. 
We assume that the commuting self-adjoint operators 
$H_0$ and $S^\ast_0S_0$ {\it strongly} commute. Then $H_0$ commutes 
with the spectral projections of $S^\ast_0 S_0$, so $H_0$ 
leaves $\Hh_{0k}$, $k\in\N_0$, invariant. Hence there are positive 
self-adjoint operators $H_{0k}$ on $\Hh_{0k}$ such that 
$H_0\zeta_k=\epsilon H_{0k} \zeta_k$, $\epsilon\in \{1,-1\}$. From 
$H_0S_0=q^2S_0H_0$ we conclude that $H_{0k}=q^{-2k} H_{00}$. Finally, the 
relation $wH_0w^\ast=q^{-1}H_0$ implies that 
$w_0H_{00}w_0^\ast=q^{-1} H_{00}$ holds 
on the Hilbert space $\Hh_{00}$. Conversely, if the latter is true, then 
the preceding formulas define operators $H_0,S_0,w$ 
fulfilling (\ref{cons22})--(\ref{cons21}). 
The representations of the relation 
$w_0H_{00}w_0^\ast= q^{-1} H_{00}$ are derived from Lemma 4.2(i). 

Now we take the {\it second series of representations} of (\ref{cons22}). 
There is a self-adjoint operator $A_{00}$ on a Hilbert space 
$\Hh_{00}$ such that $\sigma(A_{00})\sqsubseteq (q^2,1]$ and 
\begin{equation*}
S_0 \zeta_k = (1{-}q^2)^{-1} \alpha_k (A_{00}) \zeta_{k-1},~ S^\ast_0 \zeta_k = (1{-}q^2)^{-1} \alpha_{k+1} (A_{00}) \zeta_{k+1}
\end{equation*}
acting on the Hilbert space 
$\Hh_0=\displaystyle\mathop{\oplus}\nolimits^\infty_{k=-\infty} \Hh_{0k}$, where $\Hh_{0k}=\Hh_{00}$.
Arguing as in the preceding paragraph, we conclude that there exist a 
unitary operator $w_0$ and a positive self-adjoint operator $H_{00}$ on 
the Hilbert space $\Hh_{00}$ satisfying the relations 
$A_{00} H_{00}=H_{00}A_{00}$,\hfill
 $w_0 A_{00}w_0^\ast= A_{00}$ and 
$w_0 H_{00}w_0^\ast=q^{-1} H_{00}$ such that $w\zeta_k=w_0\zeta_k$ and 
$H_0\zeta_k=q^{-2k}\epsilon H_{00}\zeta_k$ for $\zeta\in\Hh_{00}$, $k\in\Z$. 
Conversely, if the latter holds, then relations 
(\ref{cons22})--(\ref{cons21}) are satisfied.

Finally we turn to the {\it third series of representions} of 
(\ref{cons22}). Then there is a unitary operator $v_0$ on $\Hh_0$ 
such that $S_0\eta =(1-q^2)^{-1} v_0\eta$, $\eta\in\Hh_0$. 
The relation $w H_0w^\ast=q^{-1} H_0$ leads to $w |H_0|w^\ast=q^{-1} |H_0|$. 
Recall that $w H_0w^\ast=q^{-1} H_0$ by (\ref{cons21}), $H_0=\epsilon |H_0|$, 
where $\epsilon\,{\in}\,\{1,-1\}$, and $\ker\,H_0{=}\{0\}$. 
By Lemma 4.2(i), there is a self-adjoint operator $H_{00}$ on a Hilbert space $\Hh_{00}$ such that $\sigma(H_{00}){\sqsubseteq} (q,1]$ and $H_0\zeta_k=q^k \epsilon H_{00}\zeta_k$, $w\zeta_k=\zeta_{k+1}$ on $\Hh_0={\displaystyle\mathop{\oplus}\nolimits^\infty_{k=-\infty}}\Hh_{0k}$, $\Hh_{0k}=\Hh_{00}$. From $H_0S_0=q^2S_0H_0$ by (\ref{cons21}) we obtain $H_0v_0=q^2v_0H_0$, hence $v_0w^{\ast 2}$ commutes with $H_0$. 
From this and the relation $v_0w=wv_0$ by (\ref{cons21}), it follows 
that there is a unitary $v_{00}$ on $\Hh_{00}$ such that 
$v_{00} H_{00}v_{00}^\ast=H_{00}$ and $v_0\zeta_k=v_{00}\zeta_{k+2},k\in\Z$. 
Conversely, the latter gives indeed a representation of relations 
(\ref{cons22})--(\ref{cons21}).

\mn
{\it Case III.} $G_0>0~ (\delta=1)$

Then relation (\ref{cons22}) reads 
\begin{equation}\label{cons3}
S^\ast_0 S_0-q^{-2}S_0S^\ast_0=q^{-2} (1-q^2)^{-1}~.
\end{equation}
Since $q<1$, relation (\ref{cons3}) has only the Fock representation by 
Lemma 4.3. Thus there is a Hilbert space $\Hh_{00}$ such that 
$\Hh_0=\displaystyle\mathop{\oplus}\nolimits^\infty_{k=0} \Hh_{0k}, \Hh_{0k} =\Hh_{00}$, and
\begin{equation*}
S_0\zeta_k=(1{-}q^2)^{-1}q^{-k-1}\lambda_{k+1} \zeta_{k+1},~ 
S^\ast_0 \zeta_k= (1{-}q^2)^{-1}q^{-k}\lambda_k\zeta_{k-1}~.
\end{equation*}
The other consistency relations (\ref{cons21}) are treated in the same 
manner as for the Fock representation in Case II. This completes the treatment 
of Step 2.

\mn
Now we bring all the above considerations together. First we insert the 
representation of the relation $w_0H_{00}w^\ast_{0}=q^{-1} H_{00}$ from 
Lemma 4.2(i). Then we put the expressions for the operators  
$S_0=H^{-1}_0 T_0, H_0$ and $w$ derived in the preceding paragraphs 
into formulas (\ref{x2})--(\ref{x1}). In doing so, we finally obtain the 
following list of $\ast$-representations of the $\ast$-algebra 
$\U_0\lti \cO(SU_q(2))$:
\begin{align*}
&(I.1)_w: \ \,&X_2 \eta_{n} &= \lambda^{-1} q^{-n}\lambda_n w\eta_{n-1},\hspace{6.5cm}\\
&       &X_0\eta_{n}  &= \lambda^{-1} q^{-n-1} \lambda_{n+1} w^\ast\eta_{n+1},\\
&       &Y_1\eta_{n}  &=0 
\ \text{ on }\, \Hh=\displaystyle\mathop{\oplus}^\infty_{n=0}\Hh_{n},\ \Hh_{n}=\K.
\end{align*}
\begin{align*}
&(I.2)_N: &X_2 \eta_{nk} &= q^{-n+k} N \eta_{n,k+1} + 
\lambda^{-1} q^{-n}\lambda_n\eta_{n-1,k+1},\hspace{3.5cm}\\
&       &X_0\eta_{nk}  &= q^{-n+k-1} N^\ast \eta_{n,k-1}+
\lambda^{-1} q^{-n-1} \lambda_{n+1} \eta_{n+1,k-1},\\
&       &Y_1\eta_{nk}  &=0 
\text{ on } \Hh=\displaystyle\mathop{\oplus}^\infty_{n=0}          
\displaystyle\mathop{\oplus}^\infty_{k=-\infty} \Hh_{nk},\ \Hh_{nk}=\K.
\end{align*}
\begin{align*}
&(II.1)_{H,\epsilon}: &X_2\eta_{nkl} &{=} 
q^{-n-2k+l+1} \lambda_k\lambda^{-1}\epsilon H\eta_{n,k-1,l}
+q^{-n}\lambda^{-1}\lambda_n\eta_{n-1,k,l+1},\\
&                   &X_0\eta_{nkl} &{=} 
q^{-n-2k+l-1} \lambda_{k+1}\lambda^{-1}\epsilon H\eta_{n,k+1,l}
+q^{-n-1} \lambda_{n+1}\lambda^{-1}\eta_{n+1,k,l-1},\\
&                   &Y_1 \eta_{nkl} &{=} {-}q^{-2n-4k+2l} H^2\eta_{nkl}
\text{ on }\Hh{=}{\displaystyle\mathop{\oplus}^\infty_{n,k=0}}~
{\displaystyle\mathop{\oplus}^\infty_{l=-\infty}} \Hh_{nkl},\ \Hh_{nkl}=\K.
\end{align*}
\begin{align*}
(II.2)_{A,H,\epsilon}\!: X_2\eta_{nkl} &{=} 
q^{-n-2k+l+1} \alpha_k(A)\lambda^{-1}\epsilon H\eta_{n,k-1,l}
+q^{-n}\lambda_n \lambda^{-1}\eta_{n-1,k,l+1},\\
                    X_0\eta_{nkl} &{=} 
q^{-n-2k+l-1} \alpha_{k+1} (A)\lambda^{-1}\epsilon H\eta_{n,k+1,l} 
{+}q^{-n-1} \lambda_{n+1} \lambda^{-1}\eta_{n+1,k,l-1},\\
                    Y_1 \eta_{nkl} &{=}{-}q^{{-}2n{-}4k{+}2l} 
H^2 \eta_{nkl}  \text{ on }\Hh{=}{\displaystyle\mathop{\oplus}^\infty_{n{=}0}}~
{\displaystyle\mathop{\oplus}^\infty_{k,l=-\infty}}\Hh_{nkl},\Hh_{nkl}=\K.
\end{align*}
\begin{align*}
&(II.3)_{H,v,\epsilon}: &X_2\eta_{nk} &{=}
q^{-n+k+1} \lambda^{-1}\epsilon Hv\eta_{n,k+2} 
+q^{-n} \lambda_n\lambda^{-1}\eta_{n-1,k+1},\hspace{2.0cm}\\
&                     & X_0\eta_{nk}&{=}
q^{-n+k-1}\lambda^{-1} \epsilon H v^\ast \eta_{n,k-2}
+q^{-n-1} \lambda_{n+1}\lambda^{-1}\eta_{n+1,k-1},\\
&                     &Y_1\eta_{nk} &{=}{-}q^{-2n+2k} H^2\eta_{nk} 
\text{ on }\Hh{=}{\displaystyle\mathop{\oplus}^\infty_{n=0}} ~
{\displaystyle\mathop{\oplus}^\infty_{k=-\infty}} \Hh_{nk},\ \Hh_{nk}=\K. 
\end{align*}
\begin{align*}
&(III)_{H,\epsilon}: &X_2\eta_{nkl} &{=} 
q^{-n+k+l} \lambda_{k+1}\lambda^{-1}\epsilon H\eta_{n,k+1,l}
+q^{-n}\lambda_n\lambda^{-1} \eta_{n-1,k,l+1},\\
&                    &X_0\eta_{nkl} &{=} 
q^{-n+k+l-1} \lambda_{k}\lambda^{-1}\epsilon H\eta_{n,k-1,l} 
+q^{-n-1} \lambda_{n+1} \lambda^{-1}\eta_{n+1,k,l-1},\\
&                    &Y_1 \eta_{nkl} &{=}q^{-2n+4k+2l} H^2 \eta_{nkl}  
\text{ on }\Hh{=}{\displaystyle\mathop{\oplus}^\infty_{n,k=0}} ~
{\displaystyle\mathop{\oplus}^\infty_{l=-\infty}} \Hh_{nkl},\ \Hh_{nkl}=\K.
\end{align*}
Here $\epsilon\in\{1,-1\}$, $N$ is a normal operator, $A$ and $H$ are 
self-adjoint operators, and $w$ and $v$ are unitaries acting on a 
Hilbert space $\K$ such that $\sigma(|N|)\sqsubseteq (q,1],
\sigma(A)\sqsubseteq(q^2,1]$ and $\sigma(H)\sqsubseteq (q,1]$. Further, 
$AH=HA$ in $(II.2)_{A,H,\epsilon}$ and $vH=Hv$ in 
$(II.3)_{H,v,\epsilon}$. 
The series $(I), (II)$ and $(III)$ correspond to the three cases 
$I, II,III$ discussed above. To complete the picture, we state 
the actions of the generators $a, b, c, d$:
\begin{align}
&(I.1)_w :\hspace{1.7cm}\notag\\
&a\eta_n{=}\lambda_n\eta_{n-1},\ 
	       d\eta_n{=}\lambda_{n+1} \eta_{n+1},\  
	       b\eta_n{=}q^{n+1} w\eta_n,\ c\eta_n{=}{-q^n}w^\ast\eta_n.
          \hspace{1.7cm}\notag\\
&(I.2)_N, (II.3)_{H,v,\epsilon}:\hspace{1.7cm}\notag\\
&a\eta_{nk}{=}\lambda_n\eta_{n-1,k},\  
d\eta_{nk}{=}\lambda_{n+1}\eta_{n+1,k},\ b\eta_{nk}=q^{n+1} \eta_{n,k+1},\ 
c\eta_{nk}{=}-q^n \eta_{n,k-1}~.\hspace{1.7cm}\notag\\
&(II.1)_{H,\epsilon},~ (II.2)_{A,H,\epsilon},~(III)_{H,\epsilon}:\hspace{1.8cm} \notag\\
\label{repu1}
&a \eta_{nkl}{=}\lambda_n\eta_{n-1,kl}, d \eta_{nkl} {=} 
\lambda_{n+1} \eta_{n+1,kl}, b\eta_{nkl}{=}q^{n{+}1}\eta_{nk,l{+}1}, 
c\eta_{nkl}{=}\makebox[0pt][l]{${-}q^n\eta_{nk,l-1}.$}
\end{align}
A representation on this list is irreducible if and only if the Hilbert space $\K$ has dimension one. In this case the parameters $N,A,H,w,v$ are complex numbers such that $|N|\in(q,1]$, $A\in(q^2,1]$, $H\in(q,1]$, $|w|=1$, and $|v|=1$. Representations corresponding to different sets of parameters $N,A,H,w,v,\epsilon$, respectively, or belonging to different series are not unitarily equivalent.

Recall that the generator $Y_1$ of $\U_0$ corresponds to the element 
$K^4$ of $\U_q(su_2)$. Hence only the representations $(III)_{H,\epsilon}$ of 
$\U_0\lti \cO(SU_q(2))$ extend to $\ast$-representations of the larger 
$\ast$-algebra $\U_q(su_2)\lti\cO(SU_q(2))$.

\subsection{Representations of the ${\bf\ast}$-algebra 
${\bf \U_q(su_2)\lti\cO(SU_q(2))}$}
The procedure is similar to that in the preceding subsection. Let us suppose that we have a $\ast$-representation of $\U_q(su_2)\lti\cO(SU_q(2))$
 on a Hilbert space $\Hh$ such that its restriction to the $\ast$-subalgebra 
$\cO(SU_q(2))$ is of the form given in 4.2.

\sn
{\bf 1. Step:} As in the preceding subsection we conclude that $\G =\{0\}$ and 
we assume that the operator K is essentially self-adjoint. 
From the relation $a K=q^{-1/2} Ka$ and $dK=q^{1/2} Kd$ it follows that 
there is an invertible self-adjoint operator $K_0$ on  $\Hh_0$ such that 
$K\eta_n=q^{-n/2}K_0\eta_n$. Since $b K=q^{-1/2} Kb$, we have 
$w K_0\eta_n=q^{-1/2} K_0w\eta_n$. The relation $dE=q^{1/2} E d+K^{-1} b$ 
implies that $E$ is of the form 
$E\eta_n=T_n\eta_n+R_n\eta_{n-1},\eta\in\D_0$, where $T_n$ and $R_n$ are 
linear operators on $\Hh_0$. From $a E=q^{-1/2} Ea$ we get 
$T_n=q^{-n/2} T_0$ and from $da E=Eda+q^{-1/2}K^{-1}ba$ we derive that 
$R_n=q^{-n/2} \lambda_n \lambda^{-1} K_0^{-1} w$. The relations 
$bE=q^{-1/2}E b$ and $KE=qEK$ yield $w T_0=q^{-1/2} T_0w$ and 
$K_0T_0=qT_0K_0$, respectively. The defining relation 
$EF-FE\equiv EE^\ast-E^\ast E=\lambda^{-1}(K^2-K^{-2})$ 
leads to $T_0T^\ast_0-T^\ast_0 T_0=\lambda^{-1} K^2_0$. 
We rewrite this in terms of $S_0:=K^{-1}_0T_0$. 

The operators $E,F,K$ act as
\begin{align}\label{opE}
E\eta_n&=q^{-n/2} K_0S_0\eta_n+q^{-n/2}\lambda_n\lambda^{-1} K_0^{-1} w\eta_{n-1}~,\\
\label{opF}
F\eta_n&=q^{-n/2} S_0^\ast K_0\eta_n +q^{-(n+1)/2} \lambda_{n+1}\lambda^{-1} 
w^\ast K^{-1}_0 \eta_{n+1}~,\\
\label{opK}
K\eta_n&=q^{-n/2} K_0\eta_n~,
\end{align}
where the operators $S_0,K_0,w$ satisfy the consistency conditions
\begin{align}\label{con2}
S^\ast_0S_0-q^{-2}S_0S^\ast_0&=(q(1-q^2))^{-1}~,\\
\label{con1}
wK_0w^\ast=q^{-1/2} K_0, K_0S_0&= qS_0K_0, w S_0w^\ast=S_0~.
\end{align}
Conversely, if (\ref{con2}) and (\ref{con1}) are fulfilled, then the operators ${E,F,K}$ defined by (\ref{opE})--(\ref{opK}) satisfy the defining relations of the algebra ${\U_q(su_2){\lti}\cO(SU_q(2))}$.

\sn
{\bf 2. Step:}
Next we investigate the consistency relations (\ref{con2}) and (\ref{con1}). 
Since $0<q<1$, by Lemma 4.3 relation (\ref{con2}) 
has only the Fock representation (see also Case III in 6.1). Hence 
there is a Hilbert space $\Hh_{00}$ such that 
$\Hh_0=\displaystyle\mathop{\oplus}\nolimits^\infty_{k=0}\Hh_{0k}, \Hh_{0k}=\Hh_{00}$, and
$$
S_0\zeta_k=q^{-1/2}\lambda^{-1}q^{-k-1}\lambda_{k+1}\zeta_{k+1},~ 
S^\ast_0\zeta_k=q^{-1/2}\lambda^{-1}q^{-k}\lambda_k\zeta_{k-1}~.
$$
Arguing as in the preceding subsection, it follows from 
(\ref{con1}) that there are a unitary operator $w_0$ and an invertible 
self-adjoint operator $K_{00}$ on $\Hh_{00}$ satisfying 
$w_0K_{00}w^\ast_0=q^{-1/2}K_{00}$ such that $w\zeta_k=w_0\zeta_k$ and  
$K_0\zeta_k=q^{k}K_{00}\zeta_k,\zeta\in\Hh_{00}$. Using the 
representations of the relation $w_0K_{00}w^\ast_0=q^{-1/2} K_{00}$ 
from Lemma 4.2(i) and inserting the preceding into the formulas 
(\ref{opE})--(\ref{opK}) we obtain the following series of 
$\ast$-representations of $\U_q(su_2)\lti\cO(SU_q(2))$:
\begin{align*}
&(I)_{H,\epsilon}:\\
&E\eta_{nkl}=q^{-(n-l+1)/2}\lambda_{k+1}\lambda^{-1} 
\epsilon H\eta_{n,k+1,l} +
q^{-(n+2k+l+1)/2}\lambda_n\lambda^{-1}\epsilon H^{-1} \eta_{n-1,k,l+1}~,\\
&F\eta_{nkl} =q^{-(n-l+1)/2} \lambda_k\lambda^{-1}
\epsilon H\eta_{n,k-1,l}+q^{-(n+2k+l+1)/2} \lambda_{n+1}\lambda^{-1}
\epsilon H^{-1}\eta_{n+1,k,l-1},
\end{align*}
$$
K\eta_{nkl}=q^{(-n+2
k+l)/2}\epsilon H\eta_{nkl}~.
$$
Here $\epsilon\in\{1,-1\}$ and $H$ is a self-adjoint operator on a 
Hilbert space $\K$ such that $\sigma (H)\sqsubseteq (q^{1/2},1]$. 
The representation Hilbert space is the direct sum 
$\Hh=\displaystyle\mathop{\oplus}^\infty_{n,k=0} 
\displaystyle\mathop{\oplus}^\infty_{l=-\infty} \Hh_{nkl}$, where 
$\Hh_{nkl}=\K$. The action of $a, b, c, d$ is 
given by (78).

The representation $(I)_{H,\epsilon}$ is irreducible if and only if $\K=\C$. Two such representation $(I)_{H,\epsilon}$ and $(I)_{H^\prime,\epsilon^\prime}$ are unitarily equivalent if and only if 
$H=H^\prime$ and $\epsilon=\epsilon^\prime$.

In the rest of this subsection we assume that $\K=\C$ and we study the 
irreducible representation $(I)_{H,\epsilon}$, $H\in(q^{1/2},1]$, more in detail. Since the representation $(I)_{H,\epsilon}$ goes into 
$(I)_{H,-\epsilon}$ if the generators $E,F,K$ of $\U_q(su_2)$ are 
replaced by their negatives, we can restrict ourselves to $(I)_{H,1}$. If 
$H\ne 1$, then the operator $K$ has eigenvalues different from 
$q^{j/2}$, $j\in\Z$, and hence the corresponding representation of 
$\U_q(su_2)$ is not integrable. 

Fix a unit vector $\eta$ of $\K=\C$. 
Let $E,F,K,K^{-1}$ denote the operators of the series $(I)_{H,1}$ defined 
on the dense domain 
$$
\D_0:=\Lin\{\eta_{nkl}\,;\,n,k\in\N_0,\ l\in\Z\}
$$ 
of $\Hh$ and let $\overline{E}, \overline{F}, \overline{K}, {\overline{K^{-1}}}$ 
denote their closures. A crucial role plays the vector 
$$
v_0:=\sum^\infty_{n=0} (-q)^n H^{2n}\eta_{n,n,-n}.
$$
{\bf Lemma 6.1} {\it For the representation 
$(I)_{H,1}$, the vector $v_0$ belongs to the intersection 
of domains $\D(\overline{E})\cap \D(\overline{F})\cap 
\D(\overline{K})\cap\D({\overline{K^{-1}}})$ and we have 
$\overline{E}v_0=0,\overline{K} v_0=Hv_0,\overline{K^{-1}} v_0=H^{-1}v_0,$
$$
\overline{F} v_0=\lambda^{-1}q^{-1/2}\sum^\infty_{n=1} (-1)^n\lambda_nH^{2n}(H-H^{-3})\eta_{n,n-1,-n}.
$$
In particular, $\overline{F} v_0=0$ if $H=1$ and 
$\overline{F} v_0\ne 0$ if $H\ne 1$.}

\mn
{\bf Proof.} We prove that $v_0\in\D(\overline{E})$ and 
$\overline{E} v_0=0$. We set
\begin{equation*}
v_{tkm}:=\sum^{k-1}_{n=0}(-1)^nq^nH^{2n}\eta_{n,n,-n}+
\sum^m_{n=0}(-1)^{k+n}t^nq^{k+n}H^{2(k+n)} \eta_{k+n, k+n,-k-n}.
\end{equation*}
for $k,m\in\N$, $t\in(0,1)$. Clearly, $v_{tkm}\in\D_0$. By the formula for $E$ we compute 
\begin{align*}
q^{1/2}\lambda E v_{tkm} =&\sum^{m-1}_{n=0} 
\left\{ ({-}1)^{k+n}(1{-}t)t^n \lambda_{k+n+1} H^{2(k+n)+1} 
\eta_{k+n,k+n+1,-k-n}\right\}\\
&+(-1)^{k+m} t^m\lambda_{k+m+1} H^{2(k+m)+1} \eta_{k+m,k+m+1,-k-m}.
\end{align*}
Using the facts that $t,H,\lambda_i$ are in $(0,1]$ we estimate
$$
{\parallel v_0 - v_{tkm}\parallel}       
\le \sum^\infty_{n=k} q^n H^{2n}+\sum^m_{n=0} t^n q^{k+n} H^{2(k+n)}  \le 2(1-q)^{-1} q^k,
$$
\begin{align*}
q\lambda^2{\parallel E v_{tkm}\parallel}^2 &
\le \sum^{m-1}_{n=0} t^{2k} 
(1-t)^2 H^{4(k+n)+2}+t^{2m} H^{4(k+m)+2}\\
&\le (1-t)(1+t)^{-1}+t^{2m}
\end{align*}
Let $\varepsilon>0$.  For large $k\in\N$ we have 
${\parallel v_0-v_{tkm}\parallel}<\varepsilon$. Now we choose 
$t\in(0,1)$ such that $(1-t)(1+t)^{-1}<\varepsilon$. Then we take 
$m\in \N$ such that $t^{2m}<\varepsilon$. Thus, we have 
$q\lambda^2{\parallel E v_{tkm}\parallel^2}<2\varepsilon$. This 
shows that $v_0\in\D(\overline{E})$ and $\overline{E}v_0=0$.
The assertion for $F$ follows in a similar manner  by a slight modification of the preceding reasoning. 
Since $K^{\pm 1} v_{tkm}=H^{\pm 1}v_{tkm}$,
we conclude that $v_0\in\D(\overline{K^{\pm 1}})$, $\overline{K} v_0=Hv_0$
 and ${\overline{K^{-1}}} v_0=H^{-1}v_0$. \hfill $\Box$

\sn
If we apply the formula for $F$ formally to the vector $\overline{F} v_0$, we obtain 
$$
\sum^\infty_{n=2} \lambda^{-2} q^{-1}
({-}1)^n\lambda_{n-1}(\lambda_n{-}\lambda_{n-1} q^2 
H^{-4})(H^2-H^{-2})(q^{-1} H^2)^n \eta_{n,n-2,-n}.
$$
If $H\ne 1$, then ${\lim\limits_{n\rightarrow\infty}} 
(\lambda_n{-}\lambda_{n-1} q^2H^{-4})=1{-}q^2H^{-4}>0$ and 
$q^{-1}H^2>1$ because $H\in(q^{1/2},1]$, so this series does not 
belong to the Hilbert space. By a more precise argument it can be 
shown that in the case $H\ne 1$ the vector $v_0$ is not in the domain of 
$\overline{F^2}$.

Now suppose that $H=1$. Then, by Lemma 6.1, 
$\overline{E} v_0=\overline{F} v_0=0$ and $\overline{K} 
v_0=\overline{K^{-1}} v_0=v_0$. These relations are the key in order 
to prove that, roughly speaking, the representation $(I)_{1,1}$ is 
the Heisenberg representation of $\U_q(su_2)\lti \cO(SU_q(2))$. In order to 
do so, we have to pass to a larger domain which contains the vector $v_0$.

\mn
{\bf Theorem 6.2} {\it There is a unique $\ast$-representation $\pi$ of 
$\U_q(su_2)\lti\cO(SU_q(2))$ on the dense domain 
$\D=\cO(SU_q(2))v_0$ of $\Hh$ such that 
$$
\pi(E)\subseteq\overline{E}, 
\pi(F)\subseteq\overline{F}, \pi(K^{\pm 1})\subseteq 
{\overline{K^{\pm 1}}}
$$ 
and 
$$
\pi(x)=x,\ \, x\in\cO(SU_q(2)),
$$ 
where all 
operators are given by the formulas of $(I)_{1,1}$. The closure of 
this representation $\pi$ is unitarily equivalent to the Heisenberg 
representation $\pi_h$ of the cross product $\ast$-algebra 
$\U_q(su_2)\lti\cO(SU_q(2))$.}

\mn 
{\bf Proof.} Recall that for $f=E$ and $x\in\cO(SU_q(2))$ relation 
(\ref{crossfx}) of the cross product algebra $\U_q(su_2)\lti\cO(SU_q(2))$ 
reads 
\begin{equation}\label{exop}
Ex=\langle K^{-1}, x_{(1)}\rangle x_{(2)}E-q^{-1}\langle 
E,x_{(1)}\rangle x_{(2)}K^{-1}.
\end{equation}
Since $(I)_{1,1}$ is a representation of $\U_q(su_2)\lti\cO(SU_q(2))$, 
this formula remains valid for the corresponding operators on $\D_0$. By the proof of Lemma 6.1, there exists a 
sequence of vectors $w_k\in\D_0$, $k\in \N$, such that 
$E w_k\rightarrow\overline{E} v_0=0$ and 
$K^{-1} w_k\rightarrow {\overline{K^{-1}}}v_0=v_0$. We apply both sides of (\ref{exop}) to $w_k$ and pass to the limit $k\rightarrow\infty$. 
Since the operators of $\cO(SU_q(2))$ are bounded, we obtain\\ 
\vspace{0.2cm}
\centerline{$
\overline{E}xv_0=-q^{-1}\langle E,x_{(1)}\rangle x_{(2)}v_0=
(x\anf S^{-1}(E))v_0.
$}
\vspace{0.1cm}
Similarly, we get\\
\vspace{0.2cm}
\centerline{$
\overline{F}xv_0=(x\anf S^{-1}(F))v_0~{\rm and}~
{\overline{K^{\pm 1}}} xv_0=(x\anf S^{-1}(K^{\pm 1}))v_0.
$}\\
Since $\anf$ is a right action of $\U_q(su_2)$ on $\cO(SU_q(2))$, it
 follows from these formulas that there is a $\ast$-representation 
 $\pi_0$ of $\U_q(su_2)$ on $\D$ such that 
 $\pi_0(f)\subseteq {\overline f}$ for $f=E,F,K,K^{-1}$ and 
\begin{equation}\label{rep}
\pi_0(f) xv_0=(x\anf S^{-1}(f))v_0,~ f\in\U_q(su_2), x\in\cO(SU_q(2)).
\end{equation}
By (\ref{rep}), there is a $\ast$-representation $\pi$ of 
$\U_q(su_2)\lti\cO(SU_q(2))$ such that $\pi(f)=\pi_0(f)$ for 
$f\in\U_q(su_2)$ and $\pi(x)=x$ for $x\in\cO(SU_q(2))$. Since 
$\pi(f)v_0=\varepsilon (f)v_0$ by (\ref{rep}), 
the linear functional $h$ on $\cO(SU_q(2))$ defined by $h(\cdot)=
{\parallel{v_0}\parallel}^{-2}\langle\pi(\cdot)v_0,v_0\rangle$  
is $\U_q(su_2)$-invariant and satisfies $h(1)=1$. By Lemma 5.2, $h$ is the Haar state of $\cO(SU_q(2))$.
By the uniqueness assertion of Proposition 5.3, 
the closure of $\pi$ is unitarily equivalent to the Heisenberg 
representation of $\U_q(su_2)\lti \cO(SU_q(2))$. Since the representation 
$(I)_{1,1}$ is irreducible, it follows from [S], Proposition 8.3.11(i), that $\D$ is dense in $\Hh$.\hfill $\Box$

\subsection{The Haar state of ${\bf \cO(SU_q(2))}$ as a partial quantum trace}
In this subsection we use the $\ast$-representation $(I)_{H,\epsilon}$ of $\U_q(su_2)\lti\cO(SU_q(2))$ to develop another approach to the Haar state $h$ of $\cO(SU_q(2))$. We fix a unit vector $\eta\in\K$ and define a linear functional $h_0$ on $\cO(SU_q(2))$ by
$$
h_0(x)=(1-q^2){\sum\nolimits^\infty_{n=0}} q^{2n}\langle x 
\eta_{n,0,-n}, \eta_{n,0,-n}\rangle,~x\in\cO(SU_q(2)).
$$
Using (\ref{repu1}) it is easy to compute $h_0(x)$ on monomials $x=a^ib^rc^s,d^jb^rc^s$ and to see that it coincides with $h(x)$ given by formula (\ref{haar}). Thus, $h_0=h$. Hence $h_0$ is $\U_q(su_2)$-invariant because $h$ is so. We will give an independent proof of the $\U_q(su_2)$-invariance of $h_0$ by using the cross product algebra $\U_q(su_2)\lti \cO(SU_q(2))$. 
From now on we assume that $\K=\C$. Recall that $H\in(q^{1/2},1]$ and $\epsilon\in\{1,-1\}$.

Let $P$ be the orthogonal projection of $\Hh$ on the closure of the subspace 
$$
\D=\Lin\{e_n:=\eta_{n,0,-n}\,;\,n\in\N_0\}
$$ 
and let $\cP$ be the set of operators $T$ defined on $\D$ for which the closure of $PTP$ is of trace class on $\Hh$. We set
$$
\Tr_P~T:=\Tr~{\overline{PTP}},\ \, T\in\cP.
$$
The main ingredient of our invariance  proof is the following partial trace property of the functional $\Tr_P$.

\mn
{\bf Lemma 6.3} {\it For $x\in\cO(SU_q(2))$ and $f\in\{EK^{-1}, FK^{-1},K^{-1}\}$, we have $fx,xf\in\cP$ and $\Tr_P~fx=\Tr_P~xf$.}

\mn
{\bf Proof.} Since $\{e_n;n\in\N_0\}$ is a orthonormal basis of $P\Hh$, we have
\begin{equation}\label{tracep}
\Tr_P~T={\sum\nolimits^\infty_{n=0}}\langle Te_n,e_n\rangle~.
\end{equation}
Set $a^{\# k}:=a^k,a^{\#(-k)}:= a^{\ast k}, c^{\# k}:=c^k, c^{\#(-k)}:=c^{\ast k}$ for $k\in\N_0$. It suffices to prove the assertion for $x=a^{\# k} c^{\# l} (c^\ast c)^j$, $k,l\in\Z,j\in\N_0$, because $\cO(SU_q(2))$ is the linear span of these elements. First let $f=K^{-1}$. If $k\ne 0$ or $l\ne 0$, then
$\Tr_PK^{-1} x=\Tr_PxK^{-1}=0$ since $K^{-1}x\eta_{n,0,-n}, xK^{-1}\eta_{n,0,-n}\in \Hh_{n-k,0,-n-l}$. If $k=l=0$, then $PK^{-1}xP e_n=PxK^{-1} P e_n=q^{(2j+1)n} e_n$. Hence the assertion holds for $f=K^{-1}$.

Let $f=EK^{-1}$. The operator $EK^{-1}$ on $\D_0$ can by written as 
$EK^{-1}=S+q^{-3/2}\lambda^{-1} c^{-1} aK^{-2}$, where $S$ maps 
$e_n=\eta_{n,0,-n}$ on multiples of $\eta_{n,1,-n}$. Hence we have $PSxP=PxSP=0$. 
Thus it suffices to prove the assertion for $f=c^{-1}aK^{-2}$. If $k\ne -1$ 
or $l\ne 1$, then $\Tr_P c^{-1} aK^{-2}x=\Tr_Pxc^{-1}aK^{-2}=0$ since 
$c^{-1} aK^{-2} x\eta_{n,0,-n}, xc^{-1} aK^{-2} \eta_{n,0,-n}\in\Hh_{n-k-1,0,-n-l+1}$. Now let $k=-1,l=1$. Then,
\begin{align*}
&c^{-1} aK^{-2} x{=}q^2(c^\ast c)^jK^{-2} {-}q^4 (c^\ast c)^{j+1}K^{-2},\\
&xc^{-1} aK^{-2}{=}q^{-2j}(c^\ast c)^j K^{-2} {-}q^{-2j}(c^\ast c)^{j+1} K^{-2}.
\end{align*}
Since $\Tr_P~(c^\ast c)^mK^{-2}=H^{-2}(1{-}q^{2m+2})^{-1}$ by (\ref{tracep}), 
it follows from these identities that $\Tr_P~ c^{-1}aK^{-2}x=\Tr_P~xc^{-1} aK^{-2}.$\\
The proof for $f=FK^{-1}$ is similar.\hfill $\Box$

\mn
{\bf Theorem 6.4}  {\it The functional $h_0(x)=(1-q^2)H^2 \Tr_PK^{-2}x$, $x\in\cO(SU_q(2))$, 
is $\U_q(su_2)$-invariant and satisfies $h_0(1)=1$. In particular, 
it is the Haar state of $\cO(SU_q(2))$.}

\mn
{\bf Proof.} Since $K^{-2} e_n=q^{2n}H^{-2} e_n$, $h_0(x)=(1{-}q^2)H^2 \Tr_P~ K^{-2} x$ by (\ref{tracep}). 
Let $x\in\cO(SU_q(2))$. There is an element $y\in\cO(SU_q(2))$ such that $xK=Ky$. Then we have
\begin{align*}
&h_0(x\anf E)=h_0(ad_R(E)x)=h_0(KxE{-}q ExK)\\
&=q^{-1}(h_0(K^2yEK^{-1})- h_0(K^2 EK^{-1} y))\\
&=q^{-1}(1{-}q^2)H^2(\Tr_PyEK^{-1}{-}\Tr_PEK^{-1}y)=0
\end{align*}
by Lemma 6.3. Similarly, $h_0(x\anf F)=0$ and $h_0(x\anf K^{\pm 1})=h_0(x)$, so $h_0$ is $\U_q(su_2)$-invariant.\hfill $\Box$

\mn
{\bf Remark 6.5} In fact, if $\{k_n\}$ and $\{l_n\}$ are arbitrary sequences from $\N_0$ and $\Z$, respectively, then the Haar state $h$ of $\cO(SU_q(2))$ can be written as
$$
h(x)=(1-q^2){\sum\nolimits^\infty_{n=0}} q^{2n} \langle x\eta_{n,k_n,l_n}, \eta_{n,k_n,l_n}\rangle.
$$

\subsection{Representations of the ${\bf \ast}$-algebra 
${\bf \U_q(gl_2)\lti\hat{\cO}(\C^2_q)}$}
We argue similarly as in the preceding two subsections and begin with
a $\ast$-representation of the $\ast$-algebra 
$\U_q(gl_2)\lti \hat{\cO}(\C^2_q)$ such that its 
restriction to $\hat{\cO}(\C^2_q)$ is admissible and hence of the form 
described in 4.2. 

Since $\G$ is the kernel of $z_1$ and $z^\ast_1$, it follows from 
(\ref{eqp1}), (\ref{eqp3}) and (\ref{eqp5}) that $\G$ is also the kernel of 
$z_2$ and $z^\ast_2$ and that the representation leaves $\G$ 
invariant. On $\G$ we have $z_1=z^\ast_1=z_2=z^\ast_2=0$ and an 
arbitrary $\ast$-representation of $\U_q(gl_2)$. Such a
$\ast$-representation of $\U_q(gl_2)\lti \hat{\cO}(\C^2_q)$ will 
be called {\it trivial}.

From now on we assume that $\G=\{0\}$. Since $wAw^\ast=A$ by assumption, 
the operator $N:=wA$ is normal. The only differences to the 
$\ast$-algebra 
 $\U_q(su_2) \lti \cO(SU_q(2))$ are the strictly positive 
operator $A$ and the additional generator $L$ which satisfies (\ref{eqp6}) 
and belongs to the center of the algebra $\U_q(gl_2)$. 
From the relation $z_2^\ast L = p L z_2^\ast$ and the fact that $L$ is 
unitary it follows that $L$ leaves each space $\Hh_n$ invariant. Hence 
there are unitaries $L_n$ on $\Hh_n$ such that $L\eta_n=L_n\eta_n$. 
From $z_1^\ast z_1L = p^2 L z_1^\ast z_1$ we get $A^2L_n= p^2 L_n A^2$ and so $L_n^\ast A L_n = p A$. Combining the latter with the relations $z_2L=pLz_2$ and $z_1L=pLz_1$ and using the fact that ${\rm ker}~ A=\{0\}$ we  derive that $L_n=L_0$ for all $n$ and that $wL_0w^\ast=L_0$. 
Since $L$ commutes with $E$, $L_0$ commutes with $T_0$ and $K_0$.
The relation $z_1E=q^{-1/2} Ez_1$ implies $AT_0=T_0A$ and $AK_0=K_0A$.
Summarizing, it follows that the operators $E,F,K,L$ act as
\begin{align}
\label{eqpe}
E\eta_n &= q^{-n/2} T_0\eta_n+
\lambda^{-1}\lambda_n q^{-n/2} K^{-1}_0 w\eta_{n-1},\\
\label{eqpf}
F\eta_n &= q^{-n/2}T^\ast_0\eta_n+
\lambda^{-1}\lambda_{n+1} q^{-(n+1)/2} w^\ast K^{-1}_0 \eta_{n+1},\\
\label{eqpk}
K\eta_n &= q^{-n/2} K_0\eta_n, L\eta_n = L_0\eta_n,
\end{align}
where $w,L_0, T_0$ and $K_0$ are operators on the Hilbert space $\Hh_0$ 
such that $w$ and $L_0$ are unitaries, $K_0$ is invertible and self-adjoint, and the conditions 
\begin{align}
\label{conseqp1}
&wK_0w^\ast=q^{-1/2} K_0,~ K_0T_0=qT_0K_0, ~wT_0w^\ast=q^{-1/2} T_0,~ 
wL_0=L_0w, \\
\label{conseqp2}
&T_0T^\ast_0-T^\ast_0T_0=\lambda^{-1}K^2_0,\\
\label{conseqp3}
&wAw^\ast=A,~ AT_0=T_0A,~ AK_0=K_0A,\\
\label{conseqp4}
&L_0 A L_0^\ast = p^{-1} A,~L_0K_0L_0^\ast = K_0,~L_0T_0L_0^\ast =T_0, \end{align}
hold. Here $A$ is a strictly positive self-adjoint operator on $\Hh_0$. Conversely, if operators $w,A, L_0, T_0$ and $K_0$ on a 
Hilbert space $\Hh_0$ are given satisfying the 
preceding conditions (say on an invariant dense domain $\D_0$ of $\Hh_0$ such that 
$K_0\D_0 =\D_0$), then above formulas define a 
$\ast$-representation of $\U_q(gl_2)\lti\hat{\cO}(\C^2_q)$. 

Note that the operators $T_0,K_0$ and 
$w$ satisfy precisely the same relations as in case of 
$\U_q(su_2)\lti\cO(SU_q(2))$. Proceeding as in 6.2 
we obtain the following list of non-trivial $\ast$-representations of 
$\U_q(gl_2)\lti\hat{\cO}(\C^2_q)$:
\begin{align*}
&(I)_{H,B,\epsilon}:\\
&E \eta_{nklj}{=}{-q^{-(n-l+1)/2}} \lambda_{k+1} \lambda^{-1}\epsilon H 
\eta_{n,k+1,l,j}{+}q^{-(n+2k+l+1)/2} \lambda_n 
\lambda^{-1}\epsilon H^{-1} \eta_{n-1,k,l+1,j},\\
&F \eta_{nklj}{=}{-q^{-(n-l+1)/2}} \lambda_k \lambda^{-1}\epsilon 
H \eta_{n,k-1,l,j}{+}q^{-(n+2k+l+1)/2} \lambda_{n+1} \lambda^{-1}\epsilon H^{-1} \eta_{n+1,k,l-1,j},\\
& K\eta_{nklj}{=}q^{-(n-2k-l)/2}\epsilon H\eta_{nklj},~~
L\eta_{nklj}{=}\eta_{nkl,j+1},\\
& z_1\eta_{nklj}{=}q^{n+1}p^j B \eta_{nk,l+1,j},~~ 
z^\ast_1 \eta_{nklj}=q^{n+1}p^j B \eta_{nk,l-1,j},\\
&z_2 \eta_{nklj}{=}\lambda_{n+1}p^j B \eta_{n+1,klj},~~
z^\ast_2 \eta_{nklj}= \lambda_n p^j B \eta_{n-1,klj},
\end{align*}
where $\epsilon \in\{-1,1\}$. The parameters $H$ and $B$ denote self-adjoint operators 
acting on a Hilbert space $\K$ such that $\sigma(H)\sqsubseteq (q^{1/2},1]$ 
and $\sigma(B)\sqsubseteq (p,1]$ if $p < 1$ resp. 
$\sigma(B)\sqsubseteq [1,p)$ if $p > 1 $. The underlying Hilbert 
space is the direct sum 
$\Hh={\displaystyle\mathop{\oplus}\limits^\infty_{n,k=0}}\,
{\displaystyle\mathop{\oplus}\limits^\infty_{l,j=-\infty}}\Hh_{nklj}$, where $\Hh_{nklj}=\K$. 
Representations corresponding to different sets $\{H,B,\epsilon\}$ of 
parameters are not unitarily equivalent. A representation of this series is 
irreducible if and only if $\K=\C$. 

The representations of the $\ast$-subalgebra 
$\U_q(su_2)\lti\hat{\cO}(\C^2_q)$ are obtained by the same formulas as above
when the last index $j$, the constants $p^j$, and the operator $L$ are omitted. In this case $B$ is 
strictly positive.

\subsection{Representations of the ${\bf \ast}$-algebra 
${\bf\U_q(gl_2)\lti \cO(\R^3_q)}$}

Suppose  we have a 
$\ast$-representation of the $\ast$-algebra 
$\U_q(gl_2)\lti \cO(\R^3_q)$.
From the defining relations of the algebra 
$\U_q(gl_2)\lti \cO(\R^3_q)$ it follows easily that 
the subspace $\K_0:= {\rm ker}~ x_2$ is invariant under all generators 
and that $x_1=x_3=0$ on $\K_0$. On $\K_0$ we can have an 
arbitrary $\ast$-representation of $\U_q(gl_2)$. Such a
$\ast$-representation of $\U_q(gl_2)\lti \cO(\R^3_q)$ is
called {\it trivial}.

\sn
{\bf 1. Step:} From now we assume that $\ker x_2 =\{0\}$ and that the restriction of the 
$\ast$-representation to $\cO(\R^3_q)$ is admissible, that is, it is of the 
form given in 4.2. 
Since $\ker x^n_1=\Hh_0+\cdots+\Hh_{n-1}$,  
the relation $x_1K=q^{-1}Kx_1$ implies that $K$ leaves the subspace 
$\Hh_0+\cdots+\Hh_{n-1}$ invariant. Since $K$ is symmetric, $K$ leaves 
$\Hh_n$ invariant so that there are operators $K_n$ on $\Hh_0$ such 
that $K\eta_n = K_n\eta_n$. The hermitian elements $K$ and $\Q_q^2$ of the
algebra $\U_q(gl_2)\lti \cO(\R^3_q)$ commute. We assume that the 
corresponding self-adjoint operators {\it strongly} commute. This implies that 
$AK_n=K_nA$ on $\Hh_0$. Combining the latter with the relation 
$x_3K=qKx_3$, we derive $K_n= q^{-n}K_0$. A slight modification of this reasoning 
shows that the operator $L$ acts as $L\eta_n=L_0\eta_n$, where $L_0$ is a unitary operator
on $\Hh_0$ such that $L_0^\ast AL_0= p^2A$. From (\ref{cr3k}) and (\ref{cr3l}) it follows that $wK_0=K_0w$ and $wL_0=L_0w$.

Using the relations $x_1E=q^{-1}Ex_1$ and $x_3E=qEx_3+q\gamma K^{-2} x_2$ 
it can be shown by induction on $n$ that $E$ maps $\Hh_n$ into 
$\Hh_{n-1}{+}\Hh_n$. Write $E\eta_n=T_n\eta_n+S_n\eta_{n-1}$, where 
$T_n$ and $S_n$ are operators on $\Hh_0$. Since $E$ commutes with 
$\Q^2_q$, we get $A^2S_n=S_nA^2$ and $A^2 T_n=T_nA^2$. We assume that $AT_n=T_nA$. From $A^2S_n=S_nA^2$ and the relation $x_1x_3E=Ex_1x_3+q^2\gamma K^{-1}x_1x_2$ we obtain by 
comparing coefficients
$$
S_n=q^{-1/2} q^{-n} \lambda^{-1}\lambda_{2n} K^{-1}_0 w.
$$
From $x_1E=q^{-1}Ex_1$ it follows that $T_n=q^{-n}T_0$. The relations 
$KE=qEK$, $LE=EL, LK=KL$ and $x_2E=Ex_2-q\gamma K^{-1} x_1$ give $K_0T_0=qT_0K_0,L_0T_0=T_0L_0, L_0K_0=K_0L_0$ and $w T_0=T_0w$, respectively. Inserting the expressions of $E$ and $K$ 
into the equation $EE^\ast-E^\ast E=\lambda^{-1}(K^2-K^{-2})$, we derive 
$T_0T^\ast_0-T^\ast_0T_0=\lambda^{-1}(K^2_0+q^{-2}K^{-2}_0)$. 

We now summarize the preceding. The operators $E,F,K$ and $L$ act as
\begin{align}\label{F0}
E \eta_n &= q^{-n} T_0 \eta_n + q^{-1/2} q^{-n} \lambda^{-1} \lambda_{2n} K^{-1}_0 w\eta_{n-1},\\
\label{F01}
F \eta_n &= q^{-n} T^\ast_0 \eta_n + q^{-3/2} q^{-n} \lambda^{-1} \lambda_{2(n+1)} K^{-1}_0 w \eta_{n+1},\\
\label{F02}
K \eta_n &= q^{-n} K_0\eta_n,~
L\eta_n=L_0 \eta_n,
\end{align}
where the operators $T_0,K_0,L_0$ and $A$ satisfy the consistency conditions
\begin{align}\label{F1}
&T^\ast_0T_0 - T_0T^\ast_0 = -\lambda^{-1} (K^2_0+q^{-2} K^{-2}_0),\\
\label{F2}
&K_0T_0=q T_0K_0,L_0T_0=T_0L_0,L_0K_0=K_0L_0,\\
\label{F3}
&AT_0=T_0A, AK_0=K_0A,L_0AL^\ast_0=p^{-2} A,\\
\label{F4}
&wT_0=T_0w,wK_0=K_0w,wL_0=L_0w, wA=Aw.
\end{align}
Conversely, if we have two self-adjoint operators $A$, $K_0$, 
a unitary operator $L_0$, a self-adjoint unitary $w$ and an operator $T_0$ on $\Hh_0$ satisfying (\ref{F1})--(\ref{F4}), then the formulas (\ref{opr3})--(\ref{opr4}) and (\ref{F0})--(\ref{F02}) define a 
$\ast$-representation of $\U_q(gl_2)\lti\cO(\R^3_q)$.

\mn
{\bf 2. Step:} Now we analyze relations (\ref{F1})--(\ref{F3}). 
Let $T_0=v|T_0|$ be the polar decomposition of the closed operator $T_0$.
From (\ref{F1}) we conclude that $\ker (T^\ast_0T_0)=\{0\}$ and so 
$\ker v=\ker |T_0|=\{0\}$. Hence $v$ is an isometry. From $K_0T_0=qT_0K_0$ we get  $K_0T_0^\ast T_0 = T_0^\ast T_0K_0$. We assume that $K_0$ and 
$T_0^\ast T_0$ are strongly commuting self-adjoint operators. Then
$K_0|T_0|=|T_0|K_0$ and hence $qvK_0 =K_0v$. Let $v = v_u\oplus v_s$ on 
$\Hh = \Hh^u \oplus \Hh^s$ be the Wold decomposition of the isometry $v$ from Lemma 4.1.

Since $v^\ast v=1$ and $T_0T_0^\ast=vT_0^\ast T_0v^\ast$, (\ref{F1}) yields
$$ 
T_0^\ast T_0v^n= vT_0^\ast T_0v^{n-1} -\lambda^{-1}(K_0^2+q^{-2}K_0^{-2})v^n,~ n \in \N~.
$$
Since $qvK_0 =K_0v$, by Lemma 4.2 $\Hh^u$ reduces $K_0$. Using this 
fact we deduce from the preceding identity that $T_0^\ast T_0$ leaves
$\Hh^u =\cap^\infty_{n=0} v^n\Hh$ invariant. We assume that $\Hh^u$ is even reducing 
for $T_0^\ast T_0$. Then $\Hh^u$ is also reducing for $T_0$ and $T_0^\ast$. 
Let us denote the restrictions of $K_0$, $T_0$ and $T_0^\ast$ by the same symbols. For the unitary part $v_u$ of $v$ it follows from (\ref{F1}) 
that 
\begin{align*}
0&<v^n_u|T_0|^2v^{\ast n}_u=
|T_0|^2+\lambda^{-1} {\sum\limits^{n-1}_{k=0}} 
(q^{-2k} K_0^2+q^{2(k-1)} K^{-2}_0)\\
&=|T_0|^2-q\lambda^{-2} 
(1-q^{2n})(q^{-2n} K^2_0+K^{-2}_0)
\end{align*} 
for all $n\in\N$. Letting $n \to \infty$ we conclude that 
the latter is only possible when $\Hh^u=\{0\}$. That is, 
we have $v =v_s$ and $\Hh=\Hh^s$. 

The operator relation $qvK_0 =K_0v$ is treated by Lemma 4.2(ii). Then 
there is an invertible self-adjoint operator $H$ on $\Hh_{00}$
such that $v\zeta_k=\zeta_{k+1}$, $K_0\zeta_k=q^kH\zeta_k$ on $\Hh_0=
{\displaystyle\mathop{\oplus}\nolimits^\infty_{k=0}} \Hh_{0k}$, $\Hh_{0k}=\Hh_{00}$.
Since $\ker T^\ast_0=\ker v^\ast=\Hh_{00}$, Equation (\ref{F1}) gives 
$|T_0|^2\zeta_0=-\lambda^{-1}(H^2+q^{-2}H^{-2})\zeta_0$. Using 
(\ref{F1}) and the fact that $K_0v=qvK_0$ we compute
\begin{align*}
|T_0|^2 \zeta_k &= |T_0|^2 v^k \zeta_0 = 
v^k (|T_0|^2 - \lambda^{-1} {\sum\limits^k_{l=1}} 
(q^{2l} K^2_0 + q^{-2(l+1)} K^{-2}_0))\zeta_0\\
&=q^{-1} \lambda^{-2} (1-q^{2(k+1)} ) (H^2_0 + 
q^{-2(k+1)} H^{-2}_0) \zeta_k~,
\end{align*}
and hence
$$
|T_0| \zeta_k = -q^{-1/2} \lambda^{-1} \lambda_{k+1} (H^2_0 + 
q^{-2(k+1)} H^{-2}_0 )^{1/2} \zeta_k~.
$$
From $L_0 T_0 = T_0 L_0, L_0 K_0 = K_0 L_0, A T_0 = T_0 A, A K_0 = 
K_0 A$ and $L_0 A L^\ast_0 = p^{-2} A$ it follows by similar 
reasoning as used in Section 6.1 that there are a unitary $L_{00}$ and a 
self-adjoint operator $A_0$ on the Hilbert space $\Hh_{00}$ satisfying 
$L_{00} H = H L_{00}, A_0 H = H A_0$ and 
$L_{00} A_0 L^\ast_{00} = p^{-2} A_0$ such that 
$L_0\zeta_k=L_{00}\zeta_k, A\zeta_k=A_0\zeta_k$.
The representations of the relation $L_{00} A_0 L^\ast_{00} = p^{-2} A_0$
are taken from Lemma 4.2(i). The relations (\ref{F4}) are treated similarly.

Carrying out the details we obtain the following list of non-trivial 
$\ast$-representations of $\U_q(gl_2)\lti \cO(\R^3_q)$:
\begin{align*}
(I)_{A,H,w}:~
&E\eta_{nkl} =-q^{-n-1/2}\lambda^{-1} \lambda_{k+1} (H^2 + q^{-2(k+1)} H^{-2})^{1/2} \eta_{n,k+1,l}\\
&\qquad~~~~~+ q^{-n-k-1/2} \lambda^{-1} \lambda_{2n} H^{-1}w \eta_{n-1,kl},\\
&F \eta_{nkl} =- q^{-n-1/2}\lambda^{-1} \lambda_k (H^2+q^{-2k} H^{-2})^{1/2} \eta_{n,k-1,l}\\
&~~~\qquad~ + q^{-n-k-3/2} \lambda^{-1} \lambda_{2(n+1)} H^{-1} w \eta_{n+1,kl},\\
&K \eta_{nkl} = q^{-n+k} H\eta_{nkl},~L \eta_{nkl} =\eta_{nk,l+1},\\
&x_1\eta_{nkl} = (1+q^2)^{-1/2} 
\lambda_{2n} p^{2l} A \eta_{n-1,kl},~
x_2\eta_{nkl} =q^{2n+1} p^{2l} Aw \eta_{nkl},\\
&x_3 \eta_{nkl} =q (1+q^2)^{-1/2} \lambda_{2(n+1)} p^{2l} 
A \eta_{n+1,kl},
\end{align*}
where $w$, $A$ and $H$ are commuting self-adjoint 
operators on a Hilbert space $\K$ such that $w$ is unitary, $H$ is invertible, 
$\sigma (A){\subseteq} (p^2, 1]$ if $p{<}1$, 
and $\sigma(A){\subseteq} (p^{-2},1]$ if $p{>}1$. The representation space is  
$\Hh=\displaystyle\mathop{\oplus}^\infty_{n,k=0} 
\displaystyle\mathop{\oplus}^\infty_{l=-\infty} \Hh_{nkl}$, 
$\Hh_{nkl}=\K$. 

If we omit the operator $L$, the last index $l$, the constants $p^{2l}$,
and assume that $A$ is a strictly self-adjoint operator, 
then the above formulas describe non-trivial $\ast$-representations of $\U_q(su_2)\lti \cO(\R^3_q)$. If  we set in addition $A=1$ and rename $x_i$ by $y_i$, then we obtain $\ast$-representations of 
$\U_q(su_2)\lti \cO(S^2_q)$. 

The cases $H=1$ and $H=q^{1/2}$ are of particular interest. Then the above representations of $\U_q(su_2)\lti \cO(\dR^3_q)$ are determined 
by the following formulas:
\begin{align*}
(I)_{A,1,w}:\quad&E\eta_{nk} = q^{-n-k-1/2}\lambda^{-1}
(-q^{-1}\lambda_{2(k+1)}\eta_{n,k+1} + \lambda_{2n} w\eta_{n-1,k}),\\
&F \eta_{nk} = q^{-n-k-1/2}\lambda^{-1}
(-\lambda_{2k}\eta_{n,k-1}+q^{-1} \lambda_{2(n+1)}w\eta_{n+1,k}),\\
&K \eta_{nk} = q^{-n+k} \eta_{nk},\\
(I)_{A,q^{1/2},w}:\quad&E\eta_{nk} =q^{-n-k-1} \lambda^{-1} (-q^{-1}  
\lambda_{k+1}\alpha_{k+2}\eta_{n,k+1} + \lambda_{n}\alpha_n w\eta_{n-1,k}),\\
&F \eta_{nk} =q^{-n-k-1}\lambda^{-1} (-\lambda_k \alpha_{k+1}  
\eta_{n,k-1}+q^{-1} \lambda_{n+1} \alpha_{n+1} w\eta_{n+1,k}),\\
&K \eta_{nk} = q^{-n+k+1/2} \eta_{nk},
\end{align*}
where $\alpha_k:= (1+q^{2k})^{1/2}$. In both cases the operators $x_1$, $x_2$ and $x_3$ act as
$$
x_1\eta_{nk} = (1+q^2)^{-1/2}\lambda_{2n}A \eta_{n-1,k},~x_2\eta_{nk} =q^{2n+1} A w \eta_{nk},$$
$$
x_3 \eta_{nk} = q(1+q^2)^{-1/2} \lambda_{2(n+1)} A \eta_{n+1,k}.
$$
Here $A$ is a strictly positive self-adjoint operator and $w$ is a self-adjoint unitary operator on a Hilbert space $\K$ such that $wAw^\ast=A$. The representation space is 
$\Hh=\displaystyle\mathop{\oplus}\limits^\infty_{n,k=0}~\Hh_{nk},~\Hh_{nk}=\K$.

In the case $\K=\dC$, $H{=}1$ the representations $(I)_{A,1,w}$  
have been constructed in [F] and [CW]. The formulas in [CW] 
can be obtained as follows. First we replace $q$ by $q^{-1}$. Then the algebra in [CW] becomes a $\ast$-subalgebra of
$\cU_q(su_2)\lti\cO(\dR^3_q)$ by setting $X^+=-x_1$, $X^-=x_3$, $X^3=x_2$, $T^+=-q^{1/2}KE$, 
$T^-=-q^{-1/2}KF$, and $T^3=\lambda^{-1}(K^4{-}1)$. The representation in [CW] is unitarily equivalent to 
$(I)_{A,1,w}$ with $A=q^{{-}1{-}2M}|z_0|$ and $w={\rm sign}\, z_0$ by the isomorphism 
$\eta_{nk}=(-1)^k|M,{-}n{+}M,k{+}n\rangle$ if $z_0>0$ and 
$\eta_{nk}=(-1)^{n{+}k}|M,{-}n{+}M,k{-}n\rangle$ if $z_0<0$.

Next we suppose that $\K{=}\dC$, $H{=}q^{1/2}$, $w{=}\pm 1$ and $A>0$. Set $\eta {:=}1$ and 
$$ v_{1/2} := {\sum\limits^\infty_{n=0}} w^n q^n \alpha_{n+1} \eta_{n,n}.
$$
Then we have $v_{-1/2}{:=} \overline{F} v_{1/2}\ne 0$, $\overline{E} v_{1/2} {=} 0$, $\overline{F}v_{-1/2} {=}0$, 
$\overline{K}v_{\pm 1/2} {=} q^{\pm 1/2} v_{\pm 1/2}$ and $\overline{K^{-1}} v_{\pm 1/2}{=}q^{\mp 1/2}v_{\pm 1/2}$. Proceeding similarly as in Theorem 6.2 it follows that $(I)_{A,q^{1/2}, w}$ leads to an irreducible closed $\ast$-representation of $\cU_q(su_2)\lti \cO(\R^3_q)$ such that its restriction to $\cU_q(su(2))$ is the direct sum of all representations $T_{l+1/2}$ with $l\in\N_0$. In particular, $(I)_{1,q^{1/2},1}$ and $(I)_{1,q^{1/2},-1}$ give inequivalent closed irreducible $\ast$-representations of $\cU_q(su_2)\lti \cO(S^2_q)$ such that their restrictions to $\cU_q(su_2)$ are integrable. Obviously, both representations are not equivalent to the Heisenberg representation.

We now close the link between the approaches in Sections 5 and 6 by describing the Heisenberg 
representation of $\U_q(su_2)\lti \cO(S^2_q)$ in terms of the representation
$(I)_{1,1,w}$. Let $\K:=\C^2$ and $\zeta:=2^{-1/2}(1,1)\in \C^2$. Let $w$ 
be the diagonal matrix with diagonal entries $1$ and $-1$. Then, 
$w^n\zeta= 2^{-1/2}(1,(-1)^n)$. Put 
$$ 
v_0:= (1-q^2)^{1/2}\sum_{n=0}^\infty q^n w^n\zeta_{n,n}.
$$  
From Theorem 6.6 below it follows that  
$h(x):= \langle x v_0,v_0 \rangle$, $x \in \cO(S^2_q)$, is the 
unique $\U_q(su_2)$-invariant state on $\cO(S^2_q)$.
 
\mn
{\bf Theorem 6.6} {\it There is a unique $\ast$-representation $\pi$ of 
$\U_q(su_2)\lti \cO(S^2_q)$ on the domain $\D=\cO(S^2_q)v_0$ such that $\pi(x)=x, x\in\cO(SU_q(2))$, and 
$$
\pi(E)\subseteq\ F^\ast, \pi(F)\subseteq\ E^\ast, \pi(K^{\pm 1})\subseteq 
({K^{\pm 1}})^\ast,
$$ 
where all operators are given by $(I)_{1,1,w}$. The closure of the representation
$\pi$ is unitarily equivalent to the Heisenberg 
representation $\pi_h$ of $\U_q(su_2)\lti\cO(S^2_q)$.}

\mn
{\bf Proof.} Since the proof is similar to the proof of Theorem 6.2, we 
sketch only 
the necessary modifications. Put $y_1^{\# n}:=y_1^n$ and $y_1^{\# -n}:=y_1^{\ast n}$ for $n \in\N$. Let $f=E,F,K,K^{-1}$. The crucial step of the proof is to show that
\begin{align}\label{yf}
\langle y_1^{\# r}y_2^s v_0,f\eta_{nk} \rangle = 
\langle ((y_1^{\# r}y_2^s) \anf S^{-1}(f^\ast))v_0, \eta_{nk} \rangle
\end{align}
for $s,n,k, \in \N_0$ and $r \in \Z$. The verification of (\ref{yf}) is straightforward. 
By (\ref{yf}), for any $x \in \cO(SU_q(2))$ 
the vector $x v_0$ is in the domain of the adjoint operator $f^\ast$ of $f$ 
and $f^\ast$ acts on $xv_0$ by
$f^\ast(xv_0)= (x \anf S^{-1}(f^\ast))v_0$. Hence it follows that there is a $\ast$-representation 
$\pi$ of $\U_q(su_2)\lti\cO(S^2_q)$ such that $\pi(x)=x$ and 
$\pi(f)(xv_0)=(x \anf S^{-1}(f))v_0$ for all $x \in \cO(SU_q(2))$ and 
$f\in \U_q(su_2)$. Since $\pi(f)v_0=\varepsilon(f) v_0$ and 
$\parallel v_0 \parallel =1$, $h$ is a 
$\U_q(su_2)$-invariant state. Hence the closure of $\pi$ is 
unitarily equivalent to the Heisenberg representation $\pi_h$ 
of $\U_q(su_2)\lti\cO(S^2_q)$. \hfill $\Box$

\sn
Let $\pi_0$ denote the $\ast$-representation $(I)_{1,1,w}$ of 
$\U_q(su_2)\lti\cO(S^2_q)$ on the domain 
$\D_0:=\Lin\{\eta_{nk};n,k\in\N_0\}$ and let $\pi_0^\ast$ be its 
adjoint representation ([S], Definition 8.1.4). Then $\pi$ is just 
the restriction of $\pi_0^\ast$ to the subdomain $\D=\cO(S^2_q)v_0$. 
Since the adjoint representation is not a $\ast$-representation in general, 
the fact that $\pi$ is $\ast$-preserving has to be proved. In fact, this was 
done by verifying formula (\ref{yf}).

{\nopagebreak\hfill \par
\end{document}